\documentclass{agtart_a}
\pdfoutput=1
\usepackage[all]{xy}

\title[Natural transformations of tensor algebras and combinatorial
groups]{Natural transformations of tensor algebras\\ and
representations of combinatorial groups } 

\author{Jelena Grbic}
\givenname{Jelena}
\surname{Grbi\'c}
\address{Department of Mathematical Sciences\\
   University of Aberdeen\\\newline
   Meston Building\\
   Aberdeen AB24 3UE\\UK}
\email{jelena@maths.abdn.ac.uk}
\urladdr{}

\author{Jie Wu}
\givenname{Jie}
\surname{Wu}
\address{Department of Mathematics\\
  National University of Singapore\\\newline
  2 Science Drive 2\\
  Singapore}
\email{matwuj@nus.edu.sg}
\urladdr{}

\volumenumber{6}
\issuenumber{}
\publicationyear{2006}
\papernumber{77}
\startpage{2189}
\endpage{2228}

\doi{}
\MR{}
\Zbl{}


\received{23 October 2006}
\revised{}
\accepted{24 October 2006}
\published{29 November 2006}
\publishedonline{29 November 2006}
\proposed{}
\seconded{}
\corresponding{}
\editor{CPR}
\version{}

\arxivreference{math.AT/0602218}

\let\xysavmatrix\xymatrix
\def\xymatrix{\disablesubscriptcorrection\xysavmatrix}
\AtBeginDocument{\let\bar\wbar\let\tilde\wtilde\let\hat\what\let\widetilde\wwtilde}
\numberwithin{equation}{section}
\def\strutt{\vrule width 0pt depth 5pt height 12pt}
\def\struttt{\vrule width 0pt depth 7pt height 14pt}
\def\co{\mskip .7mu\colon\thinspace}

\makeatletter
\def\cnewtheorem#1[#2]#3{\newtheorem{#1}{#3}[section]
\expandafter\let\csname c@#1\endcsname\c@thm}
\makeatother

\newtheorem{thm}{Theorem}[section]    
\cnewtheorem{lem}[thm]{Lemma}          
\cnewtheorem{proposition}[thm]{Proposition}
\cnewtheorem{corollary}[thm]{Corollary}
\theoremstyle{definition}
\cnewtheorem{defn}[thm]{Definition}    
\makeautorefname{defn}{Definition}
\newtheorem*{rem}{Remark}             

\newcommand{\lra}{\longrightarrow}
\newcommand{\Ker}{\mathop{\mathrm{Ker}}}
\renewcommand{\Im}{\mathop{\mathrm{Im}}}

\newcommand{\Hom}{\ensuremath{\mathrm{Hom}}}
\newcommand{\Coalg}{\ensuremath{\mathrm{Coalg}}}

\newcommand{\funct}{\mathrm{funct}}

\def\res{\mathop{\rm res}\nolimits}

\renewcommand{\H}{\ensuremath{\mathcal{H}}}

\def\Lie{\mathop{\rm Lie}\nolimits}

\newcommand{\dr}[3]{\ensuremath{#1\stackrel{#2}
{\longrightarrow}#3}}
\newcommand{\ddr}[5]{\ensuremath{#1\stackrel{#2}
{\longrightarrow}#3\stackrel{#4}{\longrightarrow}#5}}
\newcommand{\dddr}[7]{\ensuremath{#1\stackrel{#2}
{\longrightarrow}#3\stackrel{#4}{\longrightarrow}#5
\stackrel{#6}{\longrightarrow}#7}}

\def\Lie{\mathop{\rm Lie}\nolimits}
\def\id{\mathop{\rm id}\nolimits}
\def\sk{\mathop{\rm sk}\nolimits}

\let\la=\langle
\let\ra=\rangle
\def\pinch{\mathop{\rm pinch}\nolimits}
\newcommand{\calS}{\ensuremath{\mathcal{S}}}

\begin{document}

\begin{abstract}
Natural linear and coalgebra transformations of tensor algebras
are studied. The representations of certain combinatorial groups
are given. These representations are connected to natural
transformations of tensor algebras and to the groups of the
homotopy classes of maps from the James construction to loop
spaces. Applications to homotopy theory appear in a sequel \cite{GW}.
\end{abstract}

\begin{asciiabstract}
Natural linear and coalgebra transformations of tensor algebras
are studied. The representations of certain combinatorial groups
are given. These representations are connected to natural
transformations of tensor algebras and to the groups of the
homotopy classes of maps from the James construction to loop
spaces. Applications to homotopy theory appear in a sequel.
\end{asciiabstract}

\begin{htmlabstract}
Natural linear and coalgebra transformations of tensor algebras
are studied. The representations of certain combinatorial groups
are given. These representations are connected to natural
transformations of tensor algebras and to the groups of the
homotopy classes of maps from the James construction to loop
spaces. Applications to homotopy theory appear in a sequel.
\end{htmlabstract}

\maketitle

\section{Introduction}\label{sec1}
This paper has as its motivation a problem coming from classical
homotopy theory, namely, the study of the natural maps from loop
suspensions to loop spaces. The method we propose for analyzing
properties of certain unstable maps departs from classical unstable
homotopy theoretical constructions. We apply the homology functor to
natural maps from loop suspensions, and obtain certain functorial
coalgebra transformations. This approach justifies the primarily goal
of the paper, which is the study of the algebra of natural linear
transformations of tensor algebra and related groups of natural
coalgebra transformations. These algebras and groups are studied by
means of combinatorial group theory. We start with an analogue of a
non-commutative exterior algebra defined by Cohen~\cite{Cohen} and
define several new combinatorial algebras as its
generalizations. These algebras are then identified with the
corresponding algebras of natural linear transformations of tensor
algebras. On the other hand, looking at certain subgroups of the group
of units of these combinatorial algebras, we recover the Cohen group
$K_n$ and obtain its generalizations (see \fullref{sec:cohen_groups}
for the definitions). We proceed by establishing group isomorphisms
between these combinatorial groups and the corresponding groups of
natural coalgebra transformations.

Combinatorial groups, on their own, are closely related to the
motivating problem, that is, to the group of natural maps from loop
suspensions to loop spaces. This connection is explained in more
detail in the sequel to this paper \cite{GW}. There we present
solutions to some problems in classical homotopy theory concerning
natural maps from loop suspensions to loop spaces obtained by applying
the algebraic machinery developed in this paper. Related calculations
are done by exploiting the rich structure in combinatorial groups.

Before outlining the main results achieved in the paper, let us
explaining the geometrical motivation behind doing the algebra by
quoting the main results of the successive paper.

Recall that for any pointed space $X$, the James--Hopf map
$H_k\co J(X)\lra J(X^{(k)})$ is combinatorially defined by
\[
H_k(x_1x_2\ldots x_n)=\prod_{1\leq i_1<i_2<\cdots <i_k\leq
n}(x_{i_1}x_{i_2}\ldots x_{i_k})
\]
with right lexicographical order in the product. The $n$-th fold
Samelson product $\widetilde{W}_n$ on $X$ is given by the
composite
\[
\ddr{\widetilde{W}_n\co X\wedge\ldots\wedge
X}{E\wedge\ldots\wedge E}{\Omega\Sigma X\wedge\ldots\wedge
\Omega\Sigma X}{[[\, ,\, ],\ldots,]}{\Omega\Sigma X}
\]
where $E\co X\lra\Omega\Sigma X$ is the canonical inclusion and
the second map $[[\, ,\, ],\ldots,]$ is the $n$--fold commutator.
The $n$--fold Whitehead product $W_n$ on $X$ is defined as the
adjoint of the $n$--fold Samelson product $\widetilde{W}_n$.

\medskip
{\bf Theorem A}\qua {\sl Let $X$ be a pointed space with the null
homotopic reduced diagonal $\bar{\Delta}\co X\lra X\wedge X$.
Then for $n>k$,
\[
\ddr{J(X^{(n)})}{\Omega W_n}{J(X)}{H_k}{J(X^{(k)})}
\]
is null homotopic if $k$ does not divide $n$.}

\medskip
The result concerning the Barratt Conjecture can be formulated as
follows.

\medskip
{\bf Theorem B}\qua {\sl Let $X=\Sigma X'$ be a suspension and let
$f\co X\lra \Omega Y$ be a map such that $p^r[f]=0$ in the
group $[X,\Omega Y]$. Let $J(f)\co J(X)\lra \Omega Y$ be the
canonical multiplicative extension of $f$. Then the following
hold.
\begin{enumerate}
 \item The map\qquad\qquad $J(f)|_{J_n(X)}\co J_n(X)\lra \Omega Y$\\[7pt]
has order $p^{r+t}$ in $[J_n(X),\Omega Y]$ if $n<p^{t+1}$.
 \item
The composite\qquad\qquad 
$
\ddr{J_{p^{t+1}}(X)}{J(f)|_{J_{p^{t+1}}(X)}}{\Omega
Y}{p^{r+t}}{\Omega Y}
$\\[7pt]
is homotopic to the composite

\hspace*{-0cm}$\dddr{J_{p^{t+1}}(X)}{\pinch}{X^{(p^{t+1})}}
{p^{r-1}(\sum_{\tau\in\Sigma_{p^{t+1}-1}}1\wedge
\tau)}{X^{(p^{t+1})}}{\widetilde{W}_{p^{t+1}}}J(X)\stackrel{J(f)}{\lra}{\Omega
Y,}$

\vspace*{.2cm} \noindent where $p^{r+t}\co \Omega Y\lra \Omega
Y$ is the $p^{r+t}$-th power map, $\widetilde{W}_n$~is the
$n$--fold Samelson product and

\vspace*{.1cm}\hspace*{-0cm}$1\wedge \tau (x_1\wedge \cdots\wedge
x_{p^{t+1}})=x_1\wedge x_{\tau(2)}\wedge\cdots\wedge
x_{\tau(p^{t+1})}\co X^{(p^{t+1})}\lra X^{(p^{t+1})}$

\vspace*{.2cm} \noindent is the map which permutes positions.
\item Let $g=J(f)\circ \widetilde{W}_{p^{t+1}}\circ
(\sum_{\tau\in\Sigma_{p^{t+1}-1}}1\wedge\tau)\circ p^{r-1}\co
X^{(p^{t+1})}\lra\Omega Y$. Then $g$ is an equivariant map with
respect to the symmetric group action, that is,
\[
g\circ\sigma\simeq g\text{\qua for any $\sigma\in \Sigma_{p^{t+1}}$.}
\]
\end{enumerate}}

Now we are ready to outline the main results in this paper.

Through the course of the paper $R$ will be a commutative ring
with identity unless specified differently. By $\Hom_R(C,A)$ we
denote the algebra of natural linear transformations from a
coalgebra $C$ to an algebra $A$, with the multiplication given by
the convolution product. The \emph{convolution product} $f*g$ of
$f,g\co C\lra A$ is defined by
\[
\dddr{C}{\psi}{C\otimes C}{f\otimes g}{A\otimes A}{\mu}{A}
\]
where $\psi\co C\lra C\otimes C$ is the comultiplication and
$\mu\co A\otimes A\lra A$ is the multiplication.
\medskip

Let $V$ be a free $R$--module. The \emph{James (coalgebra)
filtration} $\{J_n(V)\}_{n\geq 0}$ of the tensor algebra $T(V)$ is
defined by
\begin{equation}
\label{Jamesfiltration}
 J_n(V)=\bigoplus_{j\leq n}T_j(V)
\end{equation}
for $n\geq 0$, where $T_j(V)=V^{\otimes j}$, the $j^{th}$ stage of
the tensor length filtration for $T(V)$. Let $C(V)=J_1(V)$.
\medskip

A non-commutative analogue of an exterior algebra is given by
$A^R(y_1,y_2,\negthinspace\cdots\negthinspace, y_n)$ the quotient
algebra of the tensor algebra $T(y_1,\cdots,y_n)$ over $R$ modulo
the two sided ideal generated by the monomials $y_{i_1}\cdots
y_{i_t} $ with $i_p=i_q$ for some $1\leq p<q\leq t$.
\medskip
\begin{proposition}
There is an isomorphism of algebras
\[
\theta_n\co A^R(y_1,\cdots,y_n)\lra\Hom_R(C(-)^{\otimes
n},T(-)).
\]
\end{proposition}
Furthermore, the James filtration $\{J_n(V)\}_{n\geq 0}$ induces a
cofiltration of algebras
\[
\Hom_R(T(-),T(-))\to\negthinspace\cdots\negthinspace\to
\Hom_R(J_n(-),T(-))\to\cdots\negthinspace\to \Hom_R(J_0(-),T(-))
\]
where the algebra $\Hom_R(T(-),T(-))$ is given by the inverse
limit
\[
\Hom_R(T(-),T(-))\cong\lim_n\Hom_R(J_n(-),T(-)).
\]
Let $L^R_n$ be the equalizer of the projection maps
\[
\pi_j\co A^R(y_1,\cdots,y_n)\lra A^R(y_1,\cdots,y_{n-1})
\]
for $1\leq j\leq n$. Then we have the following result.
\begin{proposition}
The map
 \[
 \theta_n=\theta_n\mid _{L_n^R}\co L^R_n\lra \Hom_R(J_n(-), T(-))
 \]
is an algebra isomorphism for $n\geq 0$.
\end{proposition}
For two $R$--modules $C$ and $D$, define their \emph{smash product}
$C\wedge D$ to be the quotient module
\[
C\wedge D=(C\otimes D)/ (C\otimes_RR\oplus R\otimes_RD).
\]
\begin{proposition}
There are combinatorial algebras $A_n^R[k]$, $^RL^{(l)}_n$ and
$^RL^{(l),(k)}_n$ (for their definitions see \fullref{sec:natural_transformations_TA}) such
that there are algebra isomorphisms
\begin{enumerate}
\item
\strutt$\Hom_R(C(-)^{\otimes n}),T(-^{\wedge k}))\cong A_n^R[k]$
\item
\strutt$\Hom_R(J_n(-^{\otimes l}),T(-))\cong\, ^RL^{(l)}_n$
\item
\strutt$\Hom_R(J_n(-^{\otimes l}),T(-^{\otimes k}))\cong\,^RL^{(l),(k)}_n$
\vspace{-6pt}
\end{enumerate}
for $1\leq n\leq\infty$.
\end{proposition}
Let $\Coalg(C,D)$ denote the group of natural coalgebra
transformations from a coalgebra $C$ to a Hopf algebra $D$ with
the multiplication given by the convolution product. The James
filtration $\{J_n(-)\}_{n\geq 0}$ induces a cofiltration of the
progroup $\Coalg(T(-),T(-))$. Recall that $C(-)=J_1(-).$

In~\cite{SW}, Selick and Wu described some properties of the
groups $\Coalg(T(-),T(\negthinspace-))$ and $\Coalg(C(-)^{\otimes
n},T(-))$. By the following theorems, we extend their results
identifying the groups of natural coalgebra transformations of the
James filtration $\{J_n(-)\}_{n\geq 0}$ and their new
generalizations with combinatorial groups introduced in
\fullref{sec1}.
\begin{thm}
\label{thm:intro4} There is an isomorphism of groups
\[
e\co K^R_n\lra \Coalg(C(-)^{\otimes n},T(-))\quad\text{for $n\geq 0$.}
\]
\end{thm}
Define $\H^R_n$ to be the equalizer of the projection
homomorphisms
\[
p_j\co K^R_n\lra K^R_{n-1}\quad \text{ for } 1\leq j\leq n.
\]

\begin{thm}
\label{thm:intro5} There is an isomorphism of groups
\[
e\co \H^R_n\lra \Coalg(J_n(-),T(-))\quad\text{for $1\leq n\leq\infty$.}
\]
\end{thm}
Having in mind the problem solved in Theorem A, we define a
generalisation $K_n^R(k)$ of $K_n$. We set an algebraic notation
which is motivated by geometry. Let
$\{x_{i_1}|x_{i_2}|\cdots|x_{i_k}\}$ be a notation for a word of
length k, in letters $x_{i_1},x_{i_2},\ldots ,x_{i_k}$. In a
successive paper these words will be related to the composite
\[
\ddr{X^n}{p_{i_1\cdots i_k}}{X^{(k)}}{E}{J(X^{(k)})}.
\]
Let $G$ be a set consisting of all the words
$\{x_{i_1}|x_{i_2}|\cdots|x_{i_k}\}$ with $1\leq i_j\leq n$ for
$1\leq j\leq k$ such that $i_s\neq i_t$ for all $1\leq s<t\leq k$.

The group $K_n^R(k)$ is defined combinatorially for any
commutative ring $R$ so that the generators are elements of $G$
and a certain set of relations, which will be discussed in detail
later on in the paper. For $k=1$, we denote $K_n^R(1)$ by $K_n^R$.
We will call $K_n^R(k)$ a \emph{Cohen group} as it is a
generalization of the combinatorial group $K_n=K_n^{\Z}(1)$
defined by Cohen~\cite{Cohen}.
\begin{thm}
\label{thm:intro6} There are combinatorial groups $K_n^R(k)$,
$\negthinspace^R\H^{(l)}_n$ and $^R\H^{(l),(k)}_n$ (see
\fullref{sec1} for their definitions) such that there are group
isomorphisms
\begin{enumerate}
\item
\strutt$\Coalg(C(-)^{\otimes n}),T(-^{\wedge k}))\cong K_n^R(k)$
\item
\strutt$\Coalg(J_n(-^{\otimes l}),T(-))\cong\,^R\H^{(l)}_n$
\item
\strutt$\Coalg(J_n(-^{\otimes l}),T(-^{\otimes k}))\cong\,^R\H^{(l),(k)}_n$
\vspace{-6pt}
\end{enumerate} for $1\leq n\leq\infty$.
\end{thm}

The compositions of the group isomorphisms in
Theorems~\ref{thm:intro4}, \ref{thm:intro5} and \ref{thm:intro6},
with the canonical inclusion of the natural coalgebra
transformations $\Coalg(C,D)$ of tensor algebras into the natural
linear transformations $\Hom_R(C,D)$ of tensor algebras give rise
to faithful representations of the introduced combinatorial groups
to the corresponding algebras of natural linear transformations of
tensor algebras.

In the sequel to this paper \cite{GW} we establish a connection between
the combinatorial groups $\H^R_n$, $^R\H_n^{(l)}$ and
$^R\H_n^{(l)(k)}$ and the groups of the homotopy classes of maps
from the topological James construction $J(X)$ of spaces $X$ with
the null homotopic reduced diagonal. We do that by first
restricting the ring $R$ to $\Z$ or $\Z/p^r$ and then constructing
injective maps:
\begin{align*}
\strutt e_X\co \H^R_n&\lra\, [J_n(X),J(X)]\\
\strutt e_X\co ^R\H^{(l)}_n&\lra\, [J_n(X^{(l)}),J(X)]\\
\strutt e_X\co ^R\H^{(l)(k)}_n&\lra\, [J_n(X^{(l)}),J(X^{(k)})]
\end{align*}
The disposition of the paper is as follows.
\fullref{sec:cohen_groups} catalogues all the combinatorial
groups of our study and states various properties they satisfy.
\fullref{sec:natural_transformations_TA} relates combinatorial
algebras to natural linear transformations of tensor algebras.
\fullref{sec:groupsandcoalgtransfof_TA} builds up to and deals
with the primary focus of the paper, that is, establishing group
isomorphisms between the combinatorial groups defined in
\fullref{sec:cohen_groups} and the corresponding groups of
functorial coalgebra transformations of tensor algebras.
\fullref{sec:representationforKn} gives a representation of
the combinatorial group $K_n^R(k)$ and relates that group to the
group of certain functorial coalgebra transformations.

\medskip
{\bf Acknowledgements}\qua The authors would like to thank Professors
Fred Cohen, John Berrick and Paul Selick for their helpful
suggestions and kind encouragement. The first author would also
like to thank Professor John Berrick and the second author for
making it possible for her to visit the National University of
Singapore for a term and providing her with such a friendly
working atmosphere.

\section{Cohen groups and their generalizations}\label{sec:cohen_groups}
In~\cite{Cohen_pub,Cohen}, Cohen defined the combinatorial group
$K^R_n(x_1,x_2,\ldots,x_n)$ for $R=\Z$ or $\Z/p^r$. Following his
approach we define the Cohen group $K^R_n(x_1,x_2,\ldots,x_n)$ for
any commutative ring $R$ with identity. Let $R_x$ denote a copy of
$R$ labeled by $x$. Write $x^r$ for the element $r\in R_x=R$, and
just $x$ for $x^1$.
\begin{defn}
The \emph{Cohen group} $K^R_n(x_1,x_2,\ldots,x_n)$ is the quotient
group of the free product $F^R_n=\coprod_{i=1}^{n} R_{x_i}$
modulo the relations
\begin{enumerate}
 \item $[[x^{r_1}_{i_1},x^{r_2}_{i_2}],x^{r_3}_{i_3}],\cdots,x^{r_l}_{i_l}]=1$ if $i_s=i_t$
 for some $1\leq s,t\leq l$;
 \item
 $[[x_{i_1}^{r_1},x_{i_2}^{r_2}],x_{i_3}^{r_3}],\cdots,x_{i_l}^{r_l}]=
 [[x_{i_1}^{r_1'},x_{i_2}^{r_2'}],x_{i_3}^{r_3'}],\cdots,x_{i_l}^{r_l'}]$
 if there is an equality\break$r_1\cdot r_2\cdot\ldots\cdot r_l=r_1'\cdot r_2'
 \cdot\ldots\cdot r_l'$.
\end{enumerate}
where $[x,y]=x^{-1}y^{-1}xy$ and $[[y_1,y_2],y_3],\cdots,y_l]$ is
an iterated commutator.
\end{defn}
For brevity of notation, we write $K^R_n$ for the Cohen group
$K^R_n(x_1,x_2,\ldots,x_n)$ when the generators are assumed, and
for $K^\Z_n$ we use Cohen's notation $K_n$. As we will show in the
sequel, the Cohen group $K_n$ is closely related to the group
$[X^n, J(X)]$ for any space $X$ such that its reduced diagonal
$\bar\Delta\co X\lra X\wedge X$ is null homotopic.

Define group homomorphisms $p_j\co F^R_n\lra F^R_{n-1}$ and
$s_j\co F^R_{n-1}\lra F^R_n$ by
\[
p_j(x^r_i)=\left\{\begin{array}{lll}
x^r_i & \text{for} & i<j\\
1 & \text{for} & i=j\\
x^r_{i-1} & \text{for} & i>j
 \end{array}\right. \quad \quad
s_j(x^r_i)=\left\{
\begin{array}{lll} x^r_i & \text{for} & i<j\\
 x^r_{i+1}& \text{for} & i\geq j
 \end{array}\right.
\]
for $1\leq j\leq n$. It follows that the composites
$\xymatrix{F^R_n\ar[r]^{p_j} & F^R_{n-1}\ar@{->>}[r]& K^R_{n-1}}$
and $\xymatrix{F^R_{n-1}\ar[r]^{s_j}& F^R_{n}\ar@{->>}[r]&
K^R_{n}}$ factor through $K^R_n$ and $K^R_{n-1}$, respectively,
inducing the homomorphisms of the Cohen groups $p_j\co
K^R_n\lra K^R_{n-1}$ and $s_j\co K^R_{n-1}\lra K^R_n$.

The next objective is to generalise the Cohen groups $K_n$, by
defining new combinatorial groups related to the group $[X^n,
J(X^{(k)})]$.
\begin{defn}
The group $K_n^R(k)$ is defined combinatorially as the quotient
group of the free product
\[
\coprod_{\begin{array}{lll} 1\leq i_j\leq n\\
1\leq j\leq k\end{array}} R_{\{x_{i_1}|x_{i_2}|\cdots|x_{i_k}\}}
\]
modulo the relations given by the following identities:
\begin{enumerate}
 \item\struttt$\{x_{i_1}|x_{i_2}|\cdots|x_{i_k}\}^r=1$\quad
if $i_s=i_t$ for some $1\leq s<t\leq k$;
\item \struttt$[[\{x_{i_1}|x_{i_2}|\cdots|x_{i_k}\}^{r_1},
\{x_{i_{k+1}}|x_{i_{k+2}}|\cdots|x_{i_{2k}}\}^{r_2},\negthinspace\cdots$\\
      \strut\hspace{2in}$\cdots\negthinspace,
\{x_{i_{(l-1)k+1}}|x_{i_{(l-1)k+2}}|\cdots|x_{i_{lk}}\}^{r_l}]=1$\\
\struttt if $i_s=i_t$ for some $1\leq s<t\leq kl$, where
$[[a_1,a_2,\cdots,a_l]=[[a_1,a_2],a_3],\cdots\negthinspace,a_l]$\\
with $[x,y]=x^{-1}y^{-1}xy$;
\item
\struttt$[[\{x_{i_1}|x_{i_2}|\cdots|x_{i_k}\}^{r_1}\negthinspace,
\{x_{i_{k+1}}|x_{i_{k+2}}|\cdots|x_{i_{2k}}\}^{r_2},\cdots\negthinspace,$\\
\strut\hspace{2in}$\cdots\negthinspace,\{x_{i_{(l-1)k+1}}|x_{i_{(l-1)k+2}}|\cdots|x_{i_{lk}}\}^{r_l}] =$\\
\struttt$[[\{x_{i_1}|x_{i_2}|\cdots|x_{i_k}\}^{r_1'},
\{x_{i_{k+1}}|x_{i_{k+2}}|\cdots|x_{i_{2k}}\}^{r_2'}\negthinspace,\cdots\negthinspace,$\\
\strut\hspace{2in}$\cdots\negthinspace,
\{x_{i_{(l-1)k+1}}|x_{i_{(l-1)k+2}}|\cdots|x_{i_{lk}}\}^{r_l'}]$\\
\struttt if $r_1\cdot r_2\cdot\ldots\cdot r_l=r_1'\cdot r_2'\cdot\ldots\cdot r_l'$.
\end{enumerate}
\end{defn}
\begin{rem}
In a similar way as for the Cohen group $K^R_n$, whenever $R=\Z$ we
will denote $K_n^R(k)$ by $K_n(k)$. It is obvious that for $k=1$,
$K_n(1)$ is the Cohen group $K_n$.
\end{rem}

Let $q$ be an integer. Notice that the group $K_n^{\Z /q}(k)$ is
the quotient group of $K_n(k)$ modulo the following additional
relations:
\[
\{x_{i_1}|x_{i_2}|\cdots|x_{i_k}\}^q=1
\]
for each generator $\{x_{i_1}|x_{i_2}|\cdots|x_{i_k}\}$.

Selick and Wu~\cite{SW} defined combinatorial groups $\H^R_n$ in
order to study the group $[J_n(X),J(X)]$ with the main goal of
obtaining information on the group $[J(X), J(X)]$. It should be
remarked that the concept (and notation) of $\H^R_n$ was invented
by Cohen. Here we recall their definition.

\begin{defn}
The group $\H^R_n$ is defined to be the equalizer of the
projections $p_j\co K^R_n(x_1,x_2,\ldots,x_n)\lra
K^R_{n-1}(x_1,x_2,\ldots,x_{n-1})$ for $1\leq j\leq n$.
\end{defn}
By definition, as $p_i\mid_{\H^R_n}=p_j\mid_{\H^R_n}$ for $1\leq
i,j\leq n$, there is a homomorphism
$d_n\co\H^R_n\lra\H^R_{n-1}$ such that the diagram
\[
\xymatrix{
\H^R_n\ar@{^(->}[r]\ar[d]^{d_n} & K^R_n\ar[d]^{p_i}\\
\H^R_{n-1}\ar@{^(->}[r] & K^R_{n-1}}
\]
commutes for each $1\leq i\leq n$.

The following lemma was proved by Selick and Wu~\cite{SW}. As the
proof is illustrative and will be used later on in the paper, we
include it here.
\begin{lem}
\label{progroupH} The homomorphism $d_n\co\H^R_n\lra\H^R_{n-1}$
is an epimorphism for each~$n$.
\end{lem}
\begin{proof}
By induction on $k$, we show that $d_{k,n}=d_{k+1}\circ\ldots\circ
d_n\co\H^R_n\lra\H^R_k$ is an epimorphism for $k\leq n$.
Clearly, $d_{1,n}$ is an epimorphism. Suppose that $d_{k-1,n}$ is
an epimorphism with $k>1$ and let $\alpha\in\H^R_k$. Since
$d_{k-1,n}\co\H^R_n\lra\H^R_{k-1}$ is onto, we assume that
$\alpha$ lies in the kernel of $d_k\co\H^R_k\lra\H^R_{k-1}$.
Let
\[
\alpha_{k,n}=\prod_{1\leq i_1<i_2<\ldots <i_{n-k}\leq n}
s_{i_{n-k}}s_{i_{n-k-1}}\ldots s_{i_1}\alpha\in K^R_n
\]
with lexicographic order from the right. Then it is routine to
check that $\alpha_{k,n}\in\H^R_n$ with
$d_{k,n}(\alpha_{k,n})=\alpha$ and hence the result.
\end{proof}
\fullref{progroupH} results in a progroup (a tower of group
epimorphisms)
\[
\xymatrix{ \H^R\ar@{->>}[r]&\ldots \ar@{->>}[r]& \H^R_n
\ar@{->>}[r]^{d_n} & \H^R_{n-1}\ar@{->>}[r]&\ldots\ar@{->>}[r]&
\H^R_1}
\]
where $\H^R$ is the group defined by the inverse limit
\[
\H^R=\lim_{d_n}\H^R_n.
\]
We recall an important description of the kernel of
$d_n\co\H^R_n\lra\H^R_{n-1}$ that will be used later on in the
paper.

Let $R$ be a commutative ring and $\bar V$ the free $R$--module
with basis $\{x_1,\ldots,x_n\}$. Then $\Lie^R(n)$ denotes the
$R$--submodule of $\bar V^{\otimes n}$ generated by the $n$--fold
commutators $[[x_{\sigma(1)},x_{\sigma(2)}],\ldots,x_{\sigma(n)}]$
for $\sigma$ a permutation in the symmetric group $\Sigma_n$ on
$n$ letters.

\begin{thm}[Cohen~\cite{Cohen}]
\label{kernel_projectionH_n} Let $\Lambda(n)$ be the kernel of the
group homomorphism $\xymatrix{d_n\co\H^R_n\ar@{->>}[r]&
\H^R_{n-1}}$. Then $\Lambda(n)$ is isomorphic to $\Lie^R(n)$ for
$R=\Z$ or $\Z/p^r$.
\end{thm}
In this paper we go a step further by defining new combinatorial
groups and algebras in order to study natural transformations of
tensor algebras. In a sequel we translate this algebraic setting
into geometry which aims at solutions of certain problems of
classical homotopy theory. With this in mind we first define two
families of groups $^R\H_n^{(l)}$ and $^R\H_n^{(l),(k)}$ that can
be seen as generalizations of $\H_n^R$ and which will shed some
light on the study of the groups $[J_n(X^{(l)}),J(X)]$ and
$[J_n(X^{(l)}),J(X^{(k)})]$.

The projection $p_j\co K^R_n(x_1,x_2,\ldots,x_n)\lra
K^R_{n-1}(x_1,x_2,\ldots,x_{n-1})$ is given by
\[
 p_j(x^r_i)=\left\{\begin{array}{lll}
x^r_i & \text{for} & i<j\\
1 & \text{for} & i=j\\
x^r_{i-1} & \text{for} & i>j
 \end{array}\right.
\]
for $1\leq j\leq l$.

Define the projective homomorphism
\begin{equation}
\label{map:proj} p_{j+\{1,\ldots,l\}}\co
K^R_{ln}(x_1,x_2,\ldots,x_{ln})\lra
K^R_{l(n-1)}(x_1,x_2,\ldots,x_{l(n-1)})
\end{equation}
by
\[
p_{j+\{1,\ldots,l\}}(x^r_i)=\left\{\begin{array}{lll}
x^r_i & \text{for} & i<j+1\\
1 & \text{for} & j+1\leq i\leq j+l\\
x^r_{i-1} & \text{for} & i>j+l
\end{array}\right.
\]
for $0\leq j\leq n-1$.
\begin{defn}
\label{H^l_n} Define $^R\H^{(l)}_n$ to be the equalizer of the
projections
\[
\xymatrix{K^R_{ln}(x_1,x_2,\ldots,x_{ln})\bigcap
\left(\bigcap_{j=0}^{n-1}\left(\bigcap^l_{i=1}\Ker
 p_{jl+i}\right)\right)\ar@<-5pt>[d]_{p_{j+\{1,\ldots,l\}}}^{..}\ar@<5pt>[d]_{.}\\
 K^R_{l(n-1)}(x_1,x_2,\ldots,x_{l(n-1)})\bigcap\left(\bigcap_{j=0}^{n-2}\left(\bigcap^l_{i=1}\Ker
 p_{jl+i}\right)\right)}
\]
for $0\leq j\leq n-1$,
\end{defn}
As the group $^R\H^{(l)}_n$ is given by the equalizer of the
projections $p_{j+\{1,\ldots,l\}}$ for $0\leq j\leq n-1$, that is,
$p_{i+\{1,\ldots,l\}}\mid_{^R\H^{(l)}_n}=p_{j+\{1,\ldots,l\}}\mid_{^R\H^{(l)}_n}$
for $0\leq i,j\leq n-1$, there is a homomorphism
$d_n\co^R\H^{(l)}_n\lra\,^R\H^{(l)}_{n-1}$ such that the
diagram
\[
\xymatrix{ ^R\H^{(l)}_n\ar@{^(->}[r]\ar[d]^{d_n} &
K^R_{ln}\ar[d]^{p_{j+\{1,\ldots,l\}}}\\
^R\H^{(l)}_{n-1}\ar@{^(->}[r] & K^R_{l(n-1)}}
\]
commutes for each $1\leq j\leq n-1$.
\begin{lem}
\label{progroupHl} The homomorphism
$d_n\co^R\H^{(l)}_n\lra\,^R\H^{(l)}_{n-1}$ is an epimorphism
for each~$n$.
\end{lem}
\begin{proof}
Noticing that the homomorphism
\[
\dr{K^R_{ln}\bigcap\left(\bigcap_{j=0}^{n-1}\left(\bigcap^l_{i=1}\Ker
 p_{jl+i}\right)\right)}{p_{j+\{1,\ldots,l\}}} {K^R_{l(n-1)}
 \bigcap\left(\bigcap_{j=0}^{n-2}\left(\bigcap^l_{i=1}\Ker
 p_{jl+i}\right)\right)}
\]
is an epimorphism for each $0\leq j\leq n-1$, the proof follows
along the lines of the proof of \fullref{progroupH}.
\end{proof}
\fullref{progroupHl} results in a progroup
\[
\xymatrix{ ^R\H^{(l)}\ar@{->>}[r]&\ldots \ar@{->>}[r] &
^R\H^{(l)}_n\ar@{->>}[r]^{d_n} & ^R\H^{(l)}_{n-1}\ar@{->>}[r] &
\ldots\ar@{->>}[r]&^R\H^{(l)}_1}
\]
where $^R\H^{(l)}$ is the group defined by the inverse limit
\[
^R\H^{(l)}=\lim_{d_n}{}^R\H^{(l)}_n.
\]
\begin{defn}
The group $^R\H_n^{(l),(k)}$ is the subgroup of $K^R_{ln}(k)$
given as the equalizer of the projections
\[
\xymatrix{
 K^R_{ln}(k)\bigcap \left(\bigcap_{j=0}^{n-1}\left(\bigcap^l_{i=1}\Ker
 p_{jl+i}\right)\right)\ar@<-5pt>[d]_{p_{j+\{1,\ldots,l\}}}^{..}\ar@<5pt>[d]_{.}\\
 K^R_{l(n-1)}(k)\bigcap \left(\bigcap_{j=0}^{n-2}\left(\bigcap^l_{i=1}\Ker
 p_{jl+i}\right)\right)}
\]
for $0\leq j\leq n-1$, where the homomorphisms $p_j$ and
$p_{j+\{1,\ldots,l\}}$ are induced on the generalized Cohen group
$K^R_{ln}(k)$ by the projections $p_{j+\{1,\ldots,l\}}$ on
$K^R_{nl}$ defined by~\eqref{map:proj}.
\end{defn}
\medskip
In an analogous fashion as for $\H^R$ and $^R\H^{(l)}$, using the
fact that $^R\H_n^{(l),(k)}$ is given as the equalizer of the
projections $p_{j+\{1,\ldots,l\}}$, there are group epimorphisms
\[
\xymatrix{d_n\co ^R\H^{(l),(k)}_n\ar@{->>}[r]&
^R\H^{(l),(k)}_{n-1}}
\]
such that the diagram
\[
\xymatrix{ ^R\H^{(l),(k)}_n\ar@{^(->}[r]\ar[d]^{d_n} & K^R_{ln}(k)
\ar[d]^{p_{j+\{1,\ldots,l\}}}\\
^R\H^{(l),(k)}_{n-1}\ar@{^(->}[r] & K^R_{l(n-1)}(k)}
\]
commutes for each $1\leq j\leq n-1$. Thus there is a progroup
\[
\xymatrix{ ^R\H^{(l),(k)}\ar@{->>}[r]&\ldots \ar@{->>}[r]&
^R\H^{(l),(k)}_n\ar@{->>}[r]^{d_n}&^R\H^{(l),(k)}_{n-1}\ar@{->>}[r]
&\ldots \ar@{->>}[r]& ^R\H^{(l),(k)}_1}
\]
where $^R\H^{(l),(k)}$ is the group defined by the inverse limit
\[
^R\H^{(l),(k)}=\lim_{d_n}\,^R\H^{(l),(k)}_n.
\]
This introduces all the combinatorial groups we will consider in
the paper.
\medskip

\section{Combinatorial algebras and natural linear\\ transformations of tensor algebras}
\label{sec:natural_transformations_TA}
\medskip

In this section the ground ring is assumed to be a commutative
ring $R$ with identity. For a Hopf algebra $H$ over $R$, denote
the comultiplication by $\psi\co H\lra H\otimes H$; the
multiplication by $\mu\co H\otimes H\lra H$; the augmentation
by $\epsilon\co H\lra R$ and the coaugmentation by $\eta\co
R\lra H$.
\medskip

Let $V$ be a free $R$--module. The \emph{James (coalgebra)
filtration} $\{J_n(V)\}_{n\geq 0}$ of the tensor algebra $T(V)$ is
defined as an $R$--module by
\[
 J_n(V)=\bigoplus_{j\leq n}T_j(V)
\]
for $n\geq 0$, where $T_j(V)$ is the $j^{th}$ stage of the tensor
word length filtration of $T(V)$. An coalgebra structure on the
filtration is given by requiring that the elements of $V$ are
primitive and then multiplicatively extend to all of $J_n(V)$.
With this coalgebra structure $J_n(V)$ is a subcoalgebra of the
primitively generated Hopf algebra $T(V)$.
\medskip

Let $C$ be a (graded) coalgebra and let $A$ be a (graded) algebra.
The \emph{convolution product} $f*g$ of $f,g\co C\lra A$ is
defined by
\[
\dddr{C}{\psi}{C\otimes C}{f\otimes g}{A\otimes A}{\mu}{A}
\]
where $\psi\co C\lra C\otimes C$ is the comultiplication and
$\mu\co A\otimes A\lra A$ is the multiplication.

Let $\Hom_R(T(\negthinspace-),T(\negthinspace-))$ denote the set
of all functorial $R$--linear maps from~$T(V)$ to itself. The
convolution product induces a multiplication in
$\Hom_R(T(\negthinspace-\negthinspace),T(\negthinspace-))$ under
which $\Hom_R(T(-),T(-))$ becomes an algebra over~$R$.
Furthermore, the James filtration
\[
J_0(V)\subseteq J_1(V)\subseteq\cdots\subseteq
J_n(V)\subseteq\cdots\subseteq T(V)
\]
induces a cofiltration of algebras
\[
\Hom_R(T(\negthinspace-),T(\negthinspace-))\rightarrow\cdots\rightarrow
\Hom_R(J_n(-),T(\negthinspace-))\rightarrow\cdots\rightarrow
\Hom_R(J_0(-),T(\negthinspace-)).
\]
The purpose of this section is on the one hand to describe
connections between combinatorial algebras and natural linear
transformations of tensor algebras, and on the other hand to
establish the context in which an already developed machinery can
be utilized to study connections between the combinatorial groups
$K^R_n$, $\H^R_n$, $^R\H_n^{(l)}$ and $^R\H_n^{(l)(k)}$ with
certain groups of functorial coalgebra transformations of $T(-)$
and $J_n(-)$. We start by recalling Cohen's
definition~\cite{Cohen} of a non-commutative analogue of the
exterior algebra on $\{ y_1,y_2,\cdots, y_n \}$.

\begin{defn}
The Cohen algebra $A^R(y_1,y_2,\cdots, y_n)$ is defined as the
quotient algebra of the tensor algebra $T(y_1,\cdots,y_n)$ over
$R$ modulo the two sided ideal generated by the monomials
$y_{i_1}\cdots y_{i_t} $ with $i_p=i_q$ for some $1\leq p<q\leq
t$.
\end{defn}

\begin{rem}
As in the case of the Cohen group $K^R_n$, Cohen~\cite{Cohen} gave
the definition of $A^R(y_1,\cdots,y_n)$ for $R=\Z$ or $\Z/p^r$,
while in our case $R$ can be any commutative ring with identity.
\end{rem}
The following result, which Cohen stated for $R=\Z$ or $\Z/p^r$
in~\cite[Theorem~1.3.2]{Cohen}, also holds for an arbitrary
commutative ring $R$.

\begin{proposition}\label{proposition2.2}
The algebra $A^R(y_1,\cdots,y_n)$  is a graded algebra over $R$
and in degree $t$, $A(y_1,\cdots,y_n)_t$ is a free $R$--module with
basis $y_{i_{\sigma(1)}}\cdots y_{i_{\sigma(t)}}$ where $\sigma\in
\Sigma_t$ acts on $\{1,2,\cdots,t\}$ and $1\leq
i_1<i_2<\cdots<i_t\leq n$.
\end{proposition}
\begin{proof}
The proof follows directly from the definition of
$A^R(y_1,\cdots,y_n)$.
\end{proof}

Recall that $C(V)$ is defined as $J_1(V)$.
\begin{defn}\label{definition4.2}
Let $V$ be a  free $R$--module. The algebra $A_n(V)$ is defined to
be the subalgebra (under convolution) of $\Hom_R(C(V)^{\otimes n},
T(V))$ generated by the elements~$y_j$ given by the composites
\[
\xymatrix{ y_j\co C(V)^{\otimes n}\ar[r]^-{p_j}&
C(V)\ar[r]^-{q}&V\ar[r]^-{E}&T(V),}
\]
for $1\leq j\leq n$, where $p_j\co C(V)^{\otimes n}\lra C(V)$
is given by
\[
p_j(x_1\otimes x_2\otimes \cdots\otimes
x_n)=\epsilon(x_1)\cdots\epsilon(x_{j-1})x_j\epsilon(x_{j+1})\cdots
\epsilon(x_n),
\]
$q\co C(V)\lra V$ is the projection and $E\co V\lra T(V)$ is
the inclusion.
\end{defn}

\begin{lem}\label{lemma4.4}
Let $V$ be a connected graded free $R$--module. Then, in the
algebra $A_n(V)$, the following equation holds for the generators
$y_j$:
\[
y_{i_1}y_{i_2}\cdots y_{i_t}=0
\]
if $i_j=i_k$ for some $j\ne k$.
\end{lem}
\begin{proof}
The lemma can be easily proved by looking at its geometrical
realisation. There exists a space $X=\vee S^{n_\alpha}$ with
$n_{\alpha}\geq1$ such that $\widetilde H_*(X;R)=V$. Thus
$H_*(X;R)\cong C(V)$ as coalgebras and $p_j\co C(V)^{\otimes
n}\lra C(V)$ is given by
\[
(p_j)_*\co H_*(X^n;R)\lra H_*(X;R),
\]
where $p_j\co X^n\lra X$ is the $j$-th projection. Thus
$y_{i_1}\cdots y_{i_t}$ is represented by the composite
\[
\ddr{H_*(X^n;R)}{(p_{i_1\cdots i_t})_*}{H_*(X^t;R)}{\gamma}
{\widetilde{H}_*(X^{(t)};R)}\xymatrix{=V^{\otimes t}\ar@{^(->}[r]
& T(V).}
\]
where $p_{i_1\cdots i_t}(x_1,\cdots,x_n)=(x_{i_1},\cdots,x_{i_t})$
and $\gamma\co X^t\lra X^{(t)}$ is the quotient map. Notice
that if $i_j=i_k$ for some~$j\ne k$, then there exists a map
$f\co X^n\lra X^{(t-1)}$ such that
\[ \bar\Delta\circ f=\gamma\circ
p_{i_1\cdots i_t}\co X^n\lra X^{(t)}
\]
where $\bar\Delta \co X^{(t-1)}\lra X^{(t)}$ is some reduced
diagonal map. Notice that
\[
(\bar\Delta)_*\co \widetilde H_*(X^{(t-1)};R)\lra \widetilde
H_*(X^{(t)};R)
\]
is zero since $X$~is a suspension. The assertion of the lemma
follows.
\end{proof}
\medskip

\begin{corollary}\label{corollary4.4}
Let $V$ be a connected graded free $R$--module. Then the map
\[
\theta_n\co A^R(y_1,\cdots,y_n)\lra A_n(V)\subseteq
\Hom_{R}(C(V)^{\otimes n},T(V))
\]
given by $\theta_n(y_j)=y_j$ is a well-defined morphism of
algebras.
\end{corollary}
\medskip
Now we show that there exists a certain connected graded free
$R$--module~$V$ such that $\theta_n\co A^R(y_1,\cdots,y_n)\lra
A_n(V)$ is a monomorphism.
\medskip

\begin{lem}\label{lemma4.6}
Let $V$ be a free $R$--module with $\dim (V)=m$. Suppose that
$m\geq n$. Then the homomorphism
\[
\theta_n\co A^R(y_1,y_2,\cdots,y_n)\lra A_n(V)
\]
is a monomorphism.
\end{lem}
\medskip
\begin{proof}
Let $X=\vee^mS^2$ be the wedge of $m$ copies of the $2$--sphere
$S^2$. Suppose that $m\geq n$. Let $V=H_2(X;R)$. A basis for the
$R$--algebra $A^R(y_1,\cdots,y_n)_k$ is given by
\[
y_{i_{\sigma(1)}}\cdots y_{i_{\sigma(k)}},
\]
where $(i_1,\cdots,i_k)$ is taken over $1\leq i_1<\cdots<i_k\leq
n$ and $\sigma$ runs over all elements in $\Sigma_k$. Notice that
$\theta_n(y_{i_{\sigma(1)}}\cdots y_{i_{\sigma(k)}})$ is
represented by the composite
\[
\ddr{H_*(X^n;R)}{p_{i_{\sigma(1)}\cdots
i_{\sigma(k)}}}{H_*(X^k;R)}{\gamma} {\widetilde
H_*(X^{(k)};R)}\xymatrix{=V^{\otimes k}\ar@{^(->}[r]& T(V).}
\]
Let $\{x_1,\cdots, x_m\}$ be a basis for $V=H_2(X;R)$. Let $1\leq
j_1<j_2<\cdots<j_k\leq n$ and let $z_1,\cdots,z_n\in C(V)$ be such
that
\begin{enumerate}
\item $z_p=1$ if $p\not\in \{j_1,\cdots,j_k\}$; \item
$z_{j_s}=x_{j_s}$.
\end{enumerate}
Then
\[
\theta_n(y_{i_{\sigma(1)}}\negthinspace\cdots
y_{i_{\sigma(k)}})(z_1\otimes z_2\otimes\cdots\otimes
z_n)\negthinspace=\negthinspace\left\{\negthinspace\negthinspace
\begin{array}{ll}
0 \negthinspace&\negthinspace\mbox{\negthinspace for\thinspace}(i_1,\negthinspace\cdots\negthinspace\negthinspace,\negthinspace i_k)\not=(j_1,\negthinspace\cdots\negthinspace, j_k)\\
x_{j_{\sigma(1)}}\cdots x_{j_{\sigma(k)}}\negthinspace&\negthinspace\mbox{\negthinspace for\thinspace}(i_1,\negthinspace\cdots\negthinspace\negthinspace,\negthinspace i_k)=(j_1,\negthinspace\cdots\negthinspace,j_k)\\
\end{array}
\right.
\]
The assertion follows.
\end{proof}

\begin{lem}[(Lemma 2.1 in \cite{SW})]\label{lemma5.2}
Let $\phi_V\co V^{\otimes n}\lra V^{\otimes m}$ be a functorial
$R$--linear map for any free $R$--module $V$ and let
$x_1,\cdots,x_n$ be $n$ homogeneous elements in $V$.
\begin{enumerate}
\item If $\dim_R(V)=n=m$, then the element
$\phi_V(x_1\otimes\cdots\otimes x_n)$ belongs to the $R$--submodule
of $V^{\otimes n}$ spanned by the elements
\[
x_{\sigma(1)}\otimes\cdots\otimes x_{\sigma(n)},
\]
where $\sigma$ runs through all elements in $\Sigma_n$. \item If
$n\not=m$, then $\phi_V$ is the zero map.
\end{enumerate}
\end{lem}

Let $\Hom_R(C(-)^{\otimes n},T(-))$ be the set of all functorial
$R$--linear maps from $C(V)^{\otimes n}$ to $T(V)$ with the
convolution product. The same object we will be sometimes denoted
by $\Hom^{\funct}_R(C(V)^{\otimes n},T(V))$.
\begin{proposition}\label{proposition4.8}
The homomorphism
\[
\theta_n\co A^R(y_1,\cdots,y_n)\lra\Hom_R(C(-)^{\otimes
n},T(-))
\]
is an isomorphism of algebras.
\end{proposition}
\begin{proof}
By \fullref{lemma4.6}, there is a free $R$--module $V$ such that
\[
\theta_n\co A^R(y_1,\cdots,y_n)\lra A_n(V)\subseteq
\Hom_R(C(V)^{\otimes n}, T(V))
\]
is a monomorphism. Therefore from the following diagram
\[
\xymatrix{
 A^R(y_1,\cdots,y_n)\ar[r]^-{\theta_n}&\Hom_R(C(-)^{\otimes n}, T(-))
 \ar[d]^{ev}\\
 & \Hom_R(C(V)^{\otimes n}, T(V))}
\]
we have that the homomorphism
\[
\theta_n\co A^R(y_1,\cdots,y_n)\lra \Hom_R(C(-)^{\otimes n},
T(-))
\]
is a monomorphism. To show that $\theta_n$ is an epimorphism
notice that for any $V$,
\begin{gather*}
C(V)^{\otimes n}=(V\oplus R)\otimes\cdots\otimes (V\oplus
R)=\bigoplus_{1\leq i_1<\cdots<i_t\leq n}V^{\otimes t}
\\
\tag*{\hbox{and}}
T(V)=\bigoplus_{m=0}^{\infty}V^{\otimes m}.
\end{gather*}
By \fullref{lemma5.2},
\[
\Hom_R(C(-)^{\otimes n},T(-))=\bigoplus_{1\leq i_1<\cdots<i_t\leq
n}\negthinspace\Hom_R^{\funct}(V^{\otimes t},V^{\otimes
t})=\bigoplus_{1\leq i_1<\cdots<i_t\leq n}\negthinspace
R(\Sigma_t).
\]
Now the assertion follows from \fullref{proposition2.2}.
\end{proof}
\medskip

\begin{defn}\label{definition5.1}
The algebra $L^R_n$ is defined to be the equalizer of the
homomorphisms
\[
\pi_j\co A^R(y_1,\cdots,y_n)\lra A^R(y_1,\cdots,y_{n-1})
\]
for $1\leq j\leq n$, where the projection map $\pi_j$ is given by
\[
\pi_j(y_k)=\left\{
\begin{array}{lll}
y_k & for & k<j\\
0 & for & k=j\\
y_{k-1}& for &k>j.
\end{array}
\right.
\]
By definition of $L^R_n$, as $\pi_i\mid_{L^R_n}=\pi_j\mid_{L^R_n}$
for $1\leq i,j\leq n$, the homomorphisms $\pi_j\co
A^R(y_1,\cdots, y_n)\lra A^R(y_1,\cdots, y_{n-1})$ induce a
homomorphism $d_n\co L^R_n\lra L^R_{n-1}$ such that the diagram
\[
\xymatrix{ L^R_{n}\ar@{^(->}[r]\ar[d]^{d_n} & A^R(y_1,\cdots,
y_{n})\ar[d]^{\pi_j}\\
L^R_{n-1}\ar@{^(->}[r]& A^R(y_1,\cdots, y_{n-1})}
\]
commutes for $1\leq j\leq n$. The algebra $L^R_{\infty}$ is
defined by the inverse limit
\[
L^R_{\infty}=\lim_{d_n}L^R_{n}.
\]
\end{defn}
\medskip

The next objective is to find a connection between the
combinatorial algebra $L^R_n$  for each $n\geq 0$ and related
algebras of natural transformations of the tensor algebra $T(-)$.
\medskip

Consider the homomorphism
\[
\theta_n\co A^R(y_1,\cdots,y_n)\lra A_n(V)\subseteq
\Hom_R(C(V)^{\otimes n},T(V)),
\]
of \fullref{lemma4.4}, where $V$ is a  free $R$--module.
Notice that $J_n(V)$ is the coequalizer of the homomorphisms
\[
\iota_j\co C(V)^{\otimes (n-1)}\lra C(V)^{\otimes n}
\]
for $1\leq j\leq n$, where $\iota_j$ is the composite
\[\xymatrix{
C(V)^{\otimes (n-1)}\ar[r]^-{\cong} & C^{\otimes (j-1)}\otimes
{R}\otimes C(V)^{\otimes (n-j)}\ar@{^(->}[r] & C^{\otimes n}.}
\]
\medskip

Let the $i$-th projection $p_i\co C(V)^{\otimes n}\lra C(V)$ be
given by
\[
p_i(x_1\otimes\cdots\otimes x_n)=
\epsilon(x_1)\cdots\epsilon(x_{i-1})x_i\epsilon(x_{i+1})\cdots\epsilon(x_n).
\]
\medskip

Then
\[
p_i\circ \iota_j=\left\{
\begin{array}{lll}
p_i&\mbox{     for}&\mbox{  $i<j$};\\
\eta_{C(V)}\circ\epsilon_{C(V)^{\otimes (n-1)}}& \mbox{
for}&\mbox{
$i=j$};\\
p_{i-1}&\mbox{   for}&\mbox{   $i>j$}.
\end{array}
\right.
\]
\medskip

Thus there is a commutative diagram
\[
\xymatrix{ A^R(y_1,\cdots, y_n)\ar[r]^-{\theta_n}\ar[d]^{\pi_j} &
\Hom_{R}(C(V)^{\otimes n}, T(V))\ar[d]^{\iota^*_j}\\
A^R(y_1,\cdots,y_{n-1})\ar[r]^-{\theta_{n-1}}&\Hom_{R}(C(V)^{\otimes
(n-1)},T(V)),}
\]
where $\pi_j$ is defined in \fullref{definition5.1}. Hence
there exists a unique homomorphism of algebras
\begin{equation}
\label{thetaonL}
 \theta_n\co L^R_n\lra\Hom_R(J_n(-),T(-))
\end{equation}
such that the diagram
\[
\xymatrix{ L^R_n\ar[r]^-{\theta_n}\ar[d] & \Hom_R(J_n(-), T(-))\ar[d]\\
A^R(y_1,\cdots,y_n)\ar[r]^-{\theta_n}\ar[d]^{\pi_j} &
\Hom_R(C(-)^{\otimes n},T(-))\ar[d]^{\iota^*_j}\\
A^R(y_1,\cdots,y_{n-1})\ar[r]^-{\theta_{n-1}}&\Hom_R(C(-)^{\otimes(n-1)},T(-))}
\]
commutes for every $1\leq j\leq n$. Furthermore, the
homomorphism~\eqref{thetaonL} preserves the projection
homomorphisms, that is, there is a commutative diagram of algebras
over $R$
\[
\xymatrix{ L^R_n\ar[r]^-{\theta_n}\ar[d]^{d_n}
&\Hom_R(J_n(-),T(-))\ar[d]\\
L^R_{n-1}\ar[r]^-{\theta_{n-1}}&\Hom_R(J_{n-1}(-),T(-)).}
\]
In the limit, this gives a homomorphism of co-filtered algebras
\[
\theta_{\infty}\co L^R_{\infty}\lra \Hom_R(T(-),T(-)).
\]

\begin{thm}\label{theoremHisomoalg}
The algebra $L^R_{\infty}$ is isomorphic to $\Hom_R(T(-),T(-))$ as
cofiltered algebras, that is, the homomorphism
\[
\theta_n\co L^R_n\lra \Hom_R(J_n(-),T(-))
\]
is an isomorphism of algebras for each $n\geq 0$.
\end{thm}

\begin{proof}
The proof is given by induction on $n$. The assertion holds
trivially in the cases $n=0$, $1$. Suppose that the assertion
holds for $n-1$ with $n\geq2$. By \fullref{lemma4.6}, the map
\[
\theta_n\co A^R(y_1,\cdots,y_n)\lra \Hom_R(C(-)^{\otimes
n},T(-))
\]
is a monomorphism. Thus the homomorphism
\[
\theta_n\co L^R_n\lra \Hom_R(J_n(-),T(-))
\]
is a monomorphism. Let $\Gamma_n$ be the kernel of the
homomorphism
\[
L^R_n\lra L^R_{n-1}
\]
and let $\tilde\Gamma_n$ be the kernel of the homomorphism
\[
\Hom_R(J_n(-),T(-))\lra \Hom_R(J_{n-1}(-),T(-)).
\]
Then
\[
\Gamma_n=\cap_{1\leq j\leq n}\ker(\pi_j) \cong
\Gamma_n(y_1,\cdots,y_n)\cong R(\Sigma_n),
\]
where $\Gamma_n(y_1,\cdots,y_n)$ is the R-submodule of
$T_n(y_1,\cdots,y_n)$ spanned by the monomials $
y_{\sigma(1)}\cdots y_{\sigma(n)} $ as $\sigma$ runs through all
elements in $\Sigma_n$, and $R(\Sigma_n)$ is a group ring. Let
$f\in \tilde\Gamma_n$, that is, $f_V\co J_n(V)\lra T(V)$ is
such that $f_V|_{J_{n-1}(V)}\co J_{n-1}(V)\lra T(V) $ is zero.
By assertion (2) of \fullref{lemma5.2}, there exists a natural
map of ~R-modules $\phi_V\co V^{\otimes n}\lra V^{\otimes n}$
such that the diagram
\[
\xymatrix{
J_n(V)\ar[r]^{f_V}\ar[d] & T(V)\\
V^{\otimes n}\ar[r]^{\phi_V} & V^{\otimes n}\ar[u]}
\]
commutes for each $V$. By assertion (1) of \fullref{lemma5.2},
we have
\[
\phi_V(x_1\cdots
x_n)=\sum_{\sigma\in\Sigma_n}k_{\sigma}x_{\sigma(1)}\cdots
x_{\sigma(n)}
\]
for some coefficients $k_{\sigma}\in R$. Let
\[
z=\sum_{\sigma\in\Sigma_n}k_{\sigma}y_{\sigma(1)}\cdots
y_{\sigma(n)}\in \gamma_n\subseteq L^R_n.
\]
Then it is routine to check that
\[
\theta_n(z)=f.
\]
Thus the map
\[
\theta_n\co \gamma_n\lra \tilde \gamma_n
\]
is an epimorphism. The assertion follows from the $5$--Lemma.
\end{proof}
\medskip

By \fullref{lemma4.6}, we have
\begin{corollary}\label{corollary3.10}
Let $V$ be a free $R$--module with $dim_RV\geq n$. Then the
homomorphism
\[
\theta_n\co L^R_n\lra \Hom_R(J_n(V),T(V))
\]
is a monomorphism.
\end{corollary}

The algebra $L^R_{\infty}$ is called the universal convolution
algebra.

\begin{defn}
The algebra $A_n^R[k]$ is defined as the quotient algebra of the
tensor algebra generated by the words
$\{y_{i_1}|y_{i_2}|\ldots|y_{i_k}\}$ with $1\leq i_j\leq n$ for
$1\leq j\leq k$ over $R$ modulo the two sided ideal generated by
the monomials
\[
\{y_{i_1}|y_{i_2}|\ldots|y_{i_k}\}\{y_{i_{k+1}}|y_{i_{k+2}}|\ldots|y_{i_{2k}}\}
\ldots \{y_{i_{(l-1)k+1}}|y_{i_{(l-1)k+2}}|\ldots|y_{i_{lk}}\}
\]
with $i_s=i_t$ for some $1\leq s<t\leq k$.
\end{defn}
For two $R$--modules $C$ and $D$, define their \emph{smash product}
$C\wedge D$ to be the quotient module
\[
C\wedge D=(C\otimes D)/ (C\otimes_RR\oplus R\otimes_RD).
\]

\begin{defn}\label{definition36}
Let $V$ be a  free $R$--module. The algebra $A^R_n[k](V)$ is
defined to be the subalgebra (under convolution) of
$\Hom_R(C(V)^{\otimes n}, T(V^{\wedge k}))$ generated by the
elements~$\{y_{i_1}|y_{i_2}|\ldots|y_{i_k}\}$ with $1\leq i_j\leq
n$ for $1\leq j\leq k$ given by the composites
\[
\dddr{\{y_{i_1}|y_{i_2}|\ldots|y_{i_k}\}\co C(V)^{\otimes n}}
{p_{i_1,\ldots,i_k}}{C(V)^{\otimes k}}{q}{V^{\wedge
k}}{E}{T(V^{\wedge k})},
\]
where $p_{i_1,\ldots,i_k}\co C(V)^{\otimes n}\lra C(V)$ is
induced by
\[
p_j(x_1\otimes x_2\otimes \cdots\otimes
x_n)=\epsilon(x_1)\cdots\epsilon(x_{j-1})x_j\epsilon(x_{j+1})\cdots
\epsilon(x_n),
\]
$q\co C(V)\lra V$ is the projection and $E\co V^{\wedge
k}\lra T(V^{\wedge k})$ is the inclusion.
\end{defn}
\medskip

\begin{lem}
\label{lemma4.4b} Let $V$ be a connected graded free $R$--module.
Then, in the algebra $A_n^R[k](V)$, the following equation holds
for the generators $\{y_{i_1}|y_{i_2}|\ldots|y_{i_k}\}$:
\[
\{y_{i_1}|y_{i_2}|\ldots|y_{i_k}\}\{y_{i_{k+1}}|y_{i_{k+2}}|\ldots|y_{i_{2k}}\}
\ldots \{y_{i_{(l-1)k+1}}|y_{i_{(l-1)k+2}}|\ldots|y_{i_{lk}}\}=0
\]
if $i_j=i_k$ for some $j\ne k$.
\end{lem}
\medskip

\begin{proof}
The proof follows along the lines of the proof of
\fullref{lemma4.4}.
\end{proof}
\medskip

\begin{corollary}\label{corollary4.4b}
Let $V$ be a connected graded free $R$--module. Then the map
\[
\theta_n\co A^R_n[k]\lra A_n^R[k](V)\subseteq
\Hom_{R}(C(V)^{\otimes n},T(V^{\wedge k}))
\]
given by
$\theta_n(\{y_{i_1}|y_{i_2}|\ldots|y_{i_k}\})=\{y_{i_1}|y_{i_2}|\ldots|y_{i_k}\}$
is a well-defined morphism of algebras.
\end{corollary}
\medskip

\begin{lem}\label{lemma4.6b}
Let $V$ be a free $R$--module with $\dim (V)=m$. Suppose that
$m\geq n$. Then the homomorphism
\[
\theta_n\co A^R_n[k]\lra A_n^R[k](V)
\]
is a monomorphism.
\end{lem}
\medskip

\begin{proof}
The proof is analogous to that of \fullref{lemma4.6}, just with
much clumsier notation.
\end{proof}
\begin{proposition}\label{proposition4.8b}
The homomorphism
\[
\theta_n\co A^R_n[k]\lra\Hom_R(C(-)^{\otimes n},T(-^{\wedge
k}))
\]
is an isomorphism of algebras.
\end{proposition}
\medskip

\begin{proof}
\hspace*{-2mm} By \fullref{lemma4.6b}, the homomorphism
 $\theta_n\co A^R_n[k]\lra\Hom_R(C(V)^{\otimes
n},$\linebreak $T(V^{\wedge k}))$ is monomorphism for a certain
choice of a free $R$--module $V$. To prove that $\theta$ is an
epimorphism, notice that
\[
C(V)^{\otimes n}=(V\oplus R)\otimes\cdots\otimes (V\oplus
R)=\bigoplus_{1\leq i_1<\cdots<i_t\leq n}(V)^{\otimes t}.
\]
Now, the assertion follows by looking at the homomorphism $\theta$
in each degree and applying \fullref{lemma5.2}.
\end{proof}
Notice that the homomorphisms
\[
\pi_{j+\{1,\ldots,l\}}\co A^R_{ln}(x_1,x_2,\ldots,x_{ln})\lra
A^R_{l(n-1)}(x_1,x_2,\ldots,x_{l(n-1)})
\]
given by
\begin{equation}
\label{projection}
\pi_{j+\{1,\ldots,l\}}(x_i)=\left\{\begin{array}{lll}
x_i & \text{for} & i<j+1\\
1 & \text{for} & j+1\leq i\leq j+l\\
x_{i-1} & \text{for} & i>j+l
\end{array}\right.
\end{equation}
for $0\leq j\leq n-1$ induce homomorphisms
\[
\pi_{j+\{1,\ldots,l\}}\co
A^R_{ln}[k](x_1,x_2,\ldots,x_{ln})\lra
A^R_{l(n-1)}[k](x_1,x_2,\ldots,x_{l(n-1)}).
\]
\begin{defn}\label{definition3.62}
The algebra $\,^RL^{(l)(k)}_n$ is defined to be the subgroup of
$A^R_{ln}[k]$ given by the equalizer of the projections
\[
\xymatrix{
 A^R_{ln}[k]\bigcap \left(\bigcap_{j=0}^{n-1}\left(\bigcap^l_{i=1}\Ker
 \pi_{jl+i}\right)\right)\ar@<-5pt>[d]_{\pi_{j+\{1,\ldots,l\}}}^{..}\ar@<5pt>[d]_{.}\\
 A^R_{l(n-1)}[k]\bigcap \left(\bigcap_{j=0}^{n-2}\left(\bigcap^l_{i=1}\Ker
 \pi_{jl+i}\right)\right)}
\]
for $0\leq j\leq n-1$.
\end{defn}
As $\,^RL^{(l)(k)}_n$ is the equalizer of the projections
$\pi_{j+\{1,\ldots,l\}}\co A^R_{ln}[k]\lra A^R_{l(n-1)}[k]$,
they induce a projective homomorphism
$d_n\co\,^RL^{(l)(k)}_n\lra\,^RL^{(l)(k)}_{n-1}$ such that the
diagram

\[
\xymatrix{ ^RL^{(l)(k)}_{n}\ar@{^(->}[r]\ar[d]^{d_n} &
A^R_{ln}[k]\ar[d]^{\pi_{j+\{1,\ldots,l\}}}\\
^RL^{(l)(k)}_{n-1}\ar@{^(->}[r]& A^R_{l(n-1)}[k]}
\]

commutes for each $1\leq j\leq n-1$. The algebra
$^RL^{(l)(k)}_{\infty}$ is defined by the inverse limit

\[
^RL^{(l)(k)}_{\infty}=\lim_{d_n}\,^RL^{(l)(k)}_{n}.
\]

The following theorem is the analogue of
\fullref{theoremHisomoalg}.
\vspace{8pt}

\begin{thm}\label{thm:Hlkisomoalg}
The algebra $ ^RL^{(l)(k)}_{\infty}$ is isomorphic to
$Hom_R(T(-^{\wedge l}\negthinspace),T(-^{\wedge k}))$ as
cofiltered algebras, that is, the homomorphism
\[
\theta_n\co\,^RL^{(l)(k)}_n\lra\Hom_R(J_n(-^{\wedge
l}),T(-^{\wedge k}))
\]
is an isomorphism of algebras for each $n\geq 0$.
\end{thm}
\vspace{8pt}

\begin{proof}
For each $n\geq 1$, there is a commutative diagram
\[
\xymatrix{ A^R_{ln}[k]\ar[r]^-{\theta_n}\ar[d]^{\pi_{j+\{1,\ldots
,l\}}} & \Hom_{R}(C(-)^{\otimes n}, T(-^{\wedge k}))
\ar[d]^{\iota^*_{j+\{1,\ldots ,l\}}}\\
A^R_{l(n-1)}[k]\ar[r]^-{\theta_{n-1}}&\Hom_{R}(C(-)^{\otimes
(n-1)},T(-^{\wedge k})),}
\]
where $\pi_{j+\{1,\ldots ,l\}}$ is induced by
projections~\eqref{projection}.
\vspace{8pt}

Hence there exists a unique
homomorphism of algebras
\[
\theta_n\co\,^RL^{(l)(k)}_n\lra\Hom_R(J_n(-^{\wedge
l}),T(-^{\wedge k}))
\]
such that the diagram
\[
\xymatrix{ ^RL^{(l)(k)}_n\ar[r]^-{\theta_n}\ar[d] &
\Hom_R(J_n(-^{\wedge l}),T(-^{\wedge k}))\ar[d]\\
A^R_{ln}[k]\bigcap\left(\bigcap_{j=0}^{n-1}\left(\bigcap^l_{i=1}\Ker
\pi_{jl+i}\right)\right) \ar[r]^-{\theta_n}
\ar[d]^{\pi_{j+\{1,\ldots,l\}}} & \Hom_R(C(-)^{\otimes
n},T(-^{\wedge k}))\ar[d]^{\iota^*_{j+\{1,\ldots,l\}}}\\
A^R_{l(n-1)}[k]\bigcap\left(\bigcap_{j=0}^{n-2}\left(\bigcap^l_{i=1}\Ker
 \pi_{jl+i}\right)\right)\ar[r]^-{\theta_{n-1}}&\Hom_R(C(-)^{\otimes(n-1)},T(-^{\wedge
 k}))}
\]
commutes for every $0\leq j\leq n-1$.

As the homomorphisms $\theta_n$ and $\theta_{n-1}$ in the bottom
two rows are isomorphisms of algebras, it follows that
$\theta_n\co\,^RL^{(l)(k)}_n\lra\Hom_R(J_n(-^{\wedge
l}),T(-^{\wedge k}))$ is an isomorphism of algebras as well.

Furthermore, there is a commutative diagram of algebras over $R$
\[
\xymatrix{ ^RL^{(l)(k)}_n\ar[r]^-{\theta_n}\ar[d]^{d_n}
&\Hom_R(J_n(-^{\wedge l}),T(-^{\wedge k}))\ar[d]\\
^RL^{(l)(k)}_{n-1}\ar[r]^-{\theta_{n-1}}&\Hom_R(J_{n-1}(-^{\wedge
l}),T(-^{\wedge k})).}
\]
This gives an isomorphism of co-filtered algebras
\[
\theta_{\infty}\co\,^RL^{(l)(k)} _{\infty}\lra
\Hom_R(T(-^{\wedge l}),T(-^{\wedge k})).
\]
The assertion follows.
\end{proof}

\begin{rem}
In the case when $k=1$ we will drop the superscript $k$ in the
above constructions. For example, $^RL^{(l)(1)}_n$ will be denoted
by $^RL^{(l)}_n$, and the algebra $\Hom_R(J_n(-^{\wedge
l}),T(-^{\wedge 1}))$ will be denoted by $\Hom_R(J_n(-^{\wedge
l}),T(-))$.
\end{rem}

\section{Combinatorial groups and natural coalgebra transformations
of tensor algebras}\label{sec:groupsandcoalgtransfof_TA}

In this section a relation between combinatorial groups and
natural coalgebra transformations of tensor algebras is
established. As a consequence, there are representations of the
combinatorial groups into algebras of natural linear
transformations of tensor algebras.

In the algebra $A^R(y_1,\cdots,y_n)$, since $y_j^2=0$, we have
$(1+ry_j)(1-ry_j)=1$ for $r\in R$ and $1\leq j\leq n$. Thus the
element $1+ry_j$ is a unit element for $r\in R$ and $1\leq j\leq
n$.

\begin{defn}
Let $R$ be a commutative ring with identity. Define the group\break
$G(A^R(y_1,\cdots,y_n))$ to be the subgroup of the group of units
in $A^R(y_1,\cdots,y_n)$ generated by the elements $1+ry_j$ for
$r\in R$ and $1\leq j\leq n$.
\end{defn}

\begin{rem}
Cohen~\cite{Cohen} proved that if $R=\Z$ or $\Z/p^r$, then
$G(A^R(y_1,\cdots,y_n))$ is isomorphic to the group $K^R_n$.
\end{rem}
Let $G$ be a group and let $x,y\in G$. We write $[x,y]$ for the
commutator $x^{-1}y^{-1}xy$. The notation $[[x_1,x_2,\cdots,x_t]$
represents the iterated commutator from left to right, that is
$[[x_1,x_2,\cdots,x_t]=[[x_1,x_2],\cdots,x_t]$.

\begin{lem}\label{lemma2.10}
Let $1\leq j_1,\cdots,j_t\leq n$ and let $r_s\in R$ for $1\leq
s\leq t$. Let $x_s=1+r_sy_{j_s}\in G(A^R(y_1,\cdots,y_n))$. Then
\begin{enumerate}
 \item $[[x_1,x_2,\cdots,x_t]=1+r_1r_2\cdots
r_t[[y_{j_1},y_{j_2},\cdots,y_{j_t}]$;
 \item Let $m=m_1\cdot m_2\cdots m_t$. Then\\
  $[[x_1^{m_1},x_2^{m_2},\cdots,x_t^{m_t}]=1+mr_1r_2\cdots
r_t[[y_{j_1},y_{j_2},\cdots,y_{j_t}]$
\end{enumerate}
where $[[x_1,x_2,\negthinspace\cdots\negthinspace,x_t]$ is the
iterated commutator in the group
$G(A^R(y_1,\negthinspace\cdots\negthinspace,y_n))$ and
$[[y_{j_1},y_{j_2},\negthinspace\cdots\negthinspace,y_{j_t}]$ is
the iterated commutator in the algebra $A^R(y_1,\cdots,y_n)$.
\end{lem}
\begin{proof}
The proof of part $(1)$ is given by induction on $t$. The
assertion holds obviously for $t\negthinspace =\negthinspace 1$.
Suppose that the assertion holds for
$t\negthinspace-1\negthinspace$. Then
$\negthinspace[[x_1,x_2,\negthinspace\cdots\negthinspace\negthinspace,x_{t-1}]$
$=1+r_1r_2\cdots r_{t-1}[[y_{j_1},y_{j_2},\cdots,y_{j_{t-1}}].$
Notice that
\[
(1+r_1r_2\cdots
r_{t-1}[[y_{j_1},y_{j_2},\cdots,y_{j_{t-1}}])(1-r_1r_2\cdots
r_{t-1}[[y_{j_1},y_{j_2},\cdots,y_{j_{t-1}}])=1.
\]
Thus $[[x_1,x_2,\cdots,x_{t-1}]^{-1}=1-r_1r_2\cdots
r_{t-1}[[y_{j_1},y_{j_2},\cdots,y_{j_{t-1}}]$ and so
\[
[[x_1,x_2,\cdots,
x_t]=[[x_1,x_2,\cdots,x_{t-1}]^{-1}x_t^{-1}[[x_1,x_2,\cdots,x_{t-1}]
x_t
\]
\[
= (1-r_1r_2\cdots r_{t-1}\alpha)(1-r_ty_{j_t})(1+r_1r_2\cdots
r_{t-1}\alpha)(1+r_ty_{j_t})
\]
\[
 =1+r_1r_2\cdots r_{t}[[y_{j_1},y_{j_2},\cdots,y_{j_{t}}],
\]
where $\alpha=[[y_{j_1},y_{j_2},\cdots,y_{j_{t-1}}]$. The
induction is finished and the assertion of part~$(1)$ follows.

The proof of part~$(2)$ is similar calculation to that in
part~$(1)$ and is done by induction on $t$.
\end{proof}

\begin{proposition}
 \label{GKiso}
For any commutative ring with identity, there is a group
isomorphism
 \begin{gather*}
 \exp: K_n^R\lra G(A^R(y_1,\cdots,y_n))
 \\
 \tag*{\hbox{given by}}
 \exp(x^r_i)=1+ry_i.
 \end{gather*}
\end{proposition}
\begin{proof}
As the proof is similar to the proof given by Cohen~\cite{Cohen}
in the case $R=\Z$ or $\Z/p^r$, we just outline the general idea.
By inspection of \fullref{lemma2.10}, the exponent map preserves
the defining relations of $K_n^R$ and thus induces a group
homomorphism. It is by definition obviously an epimorphism. That
$\exp$ is a monomorphism follows by induction on the descending
central series filtration of $K_n^R$ and noticing that $K_n^R$ is
a nilpotent group of class $n$.
\end{proof}

Let $\bar V\negthinspace=\negthinspace\la
x_1,\negthinspace\cdots\negthinspace,x_n\negthinspace\ra$ be the
free $R$--module generated by $x_1,\cdots\negthinspace,x_n$. Let
$\gamma_n^R$~be~the $R$--submodule of $\bar V^{\otimes
n}\negthinspace$ generated by the homogenous elements
$x_{\sigma(\negthinspace 1\negthinspace)}x_{\sigma(\negthinspace
2\negthinspace)}\negthinspace\cdots\negthinspace
x_{\sigma(\negthinspace n\negthinspace)}$ for $\sigma\in S_n$. Let
$\Lie^R(n)$ be the $R$--submodule of $\gamma^R_n$ generated by the
$n$--fold commutators $[[x_{\sigma(1)}, x_{\sigma(2)},\cdots
x_{\sigma(n)}]$ for $\sigma\in S_n$. Let $P_n(T(V))$ be the set of
primitive elements of $T(V)$ of tensor length $n$.

\begin{lem}\label{lemma2.11}
$\gamma^R_n\cap P_n(T(\bar V))=\Lie^R(n)$.
\end{lem}
\begin{proof}
If $R=\Z$, then $P(T(\bar V))=L(\bar V)$ is the free Lie algebra
generated by $\bar V$. Thus $\gamma^{\Z}_n\cap P_n(T(\bar
V))=\gamma^{\Z}_n\cap L_n(\bar V)= \Lie^{\Z}(n)$.

Let $p$ be a prime. If $R=\Z/p$, then $P(T(\bar V))=L^{\res}(\bar
V)$ is the free restricted Lie algebra generated by $V$. Selick
and Wu~\cite{SW} proved that $\gamma^{\Z/p}_n\cap L^{\res}_n(\bar
V)\subseteq L_n(\bar V)$. Thus
\[
\gamma^{\Z/p}_n\cap P_n(T(\bar V))=\gamma^{\Z/p}_n\cap
L^{\res}_n(\bar V)=\gamma^{\Z/p}_n\cap L_n(\bar
V)=\Lie^{\Z/p}(n).\]

Now assume $R$ is a commutative ring with identity. Let
\[
\bar\psi\co T(\bar V)\lra T(\bar V)\otimes T(\bar V)
\]
be the reduced comultiplication, that is,
$\bar\psi(x)=\psi(x)-x\otimes 1-1\otimes x$. Then $P(T(V))$ is the
kernel of $\bar\psi\co T(\bar V)\lra T(\bar V)\otimes T(\bar
V)$. Let $\psi_s\co V^{\otimes n}\lra V^{\otimes s}\otimes
V^{\otimes n-s}$ be defined by
\[
\psi(a_1a_2\cdots a_n)=\sum_{(I,J)}a_{i_1}a_{i_2}\cdots
a_{i_s}\otimes a_{j_1}a_{j_2}\cdots a_{j_{n-s}},
\]
where $I=(i_1,i_2,\cdots,i_s)$ with $i_1<i_2<\cdots<i_s$,
$J=(j_1,j_2,\cdots,j_{n-s})$ with $j_1<j_2<\cdots<j_{n-s}$ and
$(I,J)$ runs over all $(s,n-s)$--shuffles. Then
$\bar\psi|_{V^{\otimes n}}=\sum_{s=1}^{n-1}\psi_s\co V^{\otimes
n}\lra \bigoplus_{s=1}^{n-1}V^{\otimes s}\otimes V^{\otimes
n-s}\subseteq T(V)\otimes T(V)$ for any $V$ and so $\gamma^R_n\cap
P_n(T(\bar V))$ is the kernel of the map
$\bar\psi|_{\gamma^R_n}=\sum_{s=1}^{n-1}{\psi_s}|_{\gamma^R_n}\co
\gamma^R_n\lra \bigoplus_{s=1}^{n-1}\bar V^{\otimes s}\otimes \bar
V^{\otimes n-s}.$ Let $\bar V'$ be the free $\Z$--module generated
by $x_1,x_2,\cdots,x_n$. Then
\[
\bar\psi|_{\gamma^R_n}=\bar\psi|_{\gamma^{\Z}_n}\otimes_{\Z}R\co
\gamma^R_n=\gamma^{\Z}_n\otimes_{\Z}R\to\negthinspace
\bigoplus_{s=1}^{n-1}\bar V^{\otimes s}\otimes \bar V^{\otimes
n-s}\negthinspace\negthinspace=\bigoplus_{s=1}^{n-1}\bar
{V'}^{\otimes s}\otimes_{\Z}\bar {V'}^{\otimes n-s}\otimes_{\Z}R.
\]
Notice that $\Ker(\bar\psi|_{\gamma^{\Z}_n})=\Lie^{\Z}(n)$. Thus
$\bar\psi|_{\gamma^{\Z}_n}$ factors through
$\gamma^{\Z}_n/\Lie^{\Z}(n)$ and the resulting map $j\co
\gamma^{\Z}_n/\Lie^{\Z}(n)\lra \bigoplus_{s=1}^{n-1}\bar
{V'}^{\otimes s}\otimes_{\Z}\bar {V'}^{\otimes n-s}$ is a
monomorphism. Notice that
$\Ker(\bar\psi|_{\gamma^{\Z/p}_n})=\Lie^{\Z/p}(n)$ for any prime
$p$. Thus
\begin{align*}
j\otimes_{\Z}\Z/p\co
(\gamma^{\Z}_n/\Lie^{\Z}(n))\otimes_{\Z}\Z/p&=\gamma^{\Z/p}_n/\Lie^{\Z/p}(n)
\\&\lra
\bigoplus_{s=1}^{n-1}\bar {V'}^{\otimes s}\otimes_{\Z}\bar
{V'}^{\otimes n-s}\otimes_{\Z}\Z/p
\end{align*}
is a monomorphism for any prime $p$ and so the cokernel of
\[
j\co \gamma^{\Z}_n/\Lie^{\Z}(n)\lra \bigoplus_{s=1}^{n-1}\bar
{V'}^{\otimes s}\otimes_{\Z}\bar {V'}^{\otimes n-s}
\]
is a torsion free finitely generated $\Z$--module. Thus
\[
j\otimes_{\Z}R\co (\gamma^{\Z}_n/\Lie^{\Z}(n))\otimes_{\Z}R\lra
\bigoplus_{s=1}^{n-1}\bar {V'}^{\otimes s}\otimes_{\Z}\bar
{V'}^{\otimes n-s}\otimes_{\Z}R
\]
is a monomorphism for any ring $R$ and so the image of
\[
\bar\psi|_{\gamma^R_n}\co \gamma^R_n\lra
\bigoplus_{s=1}^{n-1}\bar {V}^{\otimes s}\otimes_R\bar
{V}^{\otimes n-s}
\]
is isomorphic to $(\gamma^{\Z}_n/\Lie^{\Z}(n))\otimes_{\Z}R$.
Notice that $\gamma^{\Z}_n/\Lie^{\Z}(n)$ is a free $\Z$--module.
From the short exact sequence
\[
\xymatrix{ \Lie^{\Z}(n)\ar@{^(->}[r]&\gamma^{\Z}_n\ar@{->>}[r]
&\gamma^{\Z}_n/\Lie^{\Z}(n),}
\]
we get the short exact sequence
\[
\xymatrix{\Lie^R(n)=\Lie^{\Z}(n)\otimes_{\Z}R\ar@{^(->}[r] &
\gamma^R_n=\gamma^{\Z}_n\otimes_{\Z}R \ar@{->>}[r] &
\gamma^{\Z}_n/\Lie^{\Z}(n)\otimes_{\Z}R.}
\]
The assertion follows.
\end{proof}

\subsection*{The group of natural coalgebra transformations
of $\{J_n(-)\}_{n\geq 0}$}

Now we determine the set of functorial coalgebra maps from
$C(V)^{\otimes n}$ to $T(V)$. Let $\Coalg(C(-)^{\otimes n},T(-))$
be the set of all functorial maps of coalgebras from
$C(V)^{\otimes n}$ to $T(V)$, with the multiplication given by the
convolution product. Notice that $\Coalg(C(-)^{\otimes n},T(-))$
is a group with the convolution product.

\begin{thm}\label{theorem2.12}
There is an isomorphism of groups
\[
e\co K_n^R\lra\Coalg(C(-)^{\otimes n},T(-))\quad\text{for each $n\geq0$.}
\]
\end{thm}
\begin{proof}
Let $\theta_n\co A^R(y_1,\cdots,y_n)\lra \Hom_R(C(-)^{\otimes
n},T(-))$ be the algebra homomorphism in
\fullref{proposition4.8}. Notice that
\[
\theta_n(1+ry_j)\co C(V)^{\otimes n}\lra T(V)
\]
is a functorial map of coalgebras for $r\in R$ and $1\leq j\leq
n$. Thus
\[
\theta_n(G(A^R(y_1,\cdots,y_n)))\subseteq \Coalg(C(-)^{\otimes
n},T(-))
\]
and so
\[
e=\theta_n|_{G(A^R(y_1,\cdots,y_n))}\co
G(A^R(y_1,\cdots,y_n))\lra\Coalg(C(-)^{\otimes n},T(-))
\]
is a well-defined monomorphism of groups. We need to show that $e$
is an epimorphism. By \fullref{GKiso}, we can identify
$K_n^R$ with $G(A^R(y_1,\cdots,y_n))$.

By assuming that $V$ is a connected module, $C(V)$ is a connected
graded coalgebra. Let $\sk_t C(V)^{\otimes n}=\bigoplus_{j\leq
2t}(C(V)^{\otimes n})_j$. Then $\{\sk_t C(V)^{\otimes n}\}$ is a
(finite) coalgebra filtration of $C(V)^{\otimes n}$. Let
$\{\Gamma^tK_n^R\}$ be the descending central series of $K_n^R$
starting with $\Gamma^1K_n^R=K_n^R.$ Notice that by
\fullref{lemma2.10},
\[
\Gamma^tK_n^R\subseteq \{1\}+I^tA^R(y_1,y_2,\cdots,y_n)
\]
for each $t\geq 1$, where $I^tA^R(y_1,y_2,\cdots,y_n)$ is the
$t$--fold product of the augmentation ideal
$IA^R(y_1,y_2,\cdots,y_n)$. Let $\sk_t$ denote $\sk_tC(-)^{\otimes
n}$. Then the composite
\[
\xymatrix{
\Gamma^sK_n^R\ar@{^(->}[r]&K_n^R\ar[r]^-{e}&\Coalg(C(-)^{\otimes
n},T(-))\ar[r]&\Coalg(\sk_t,T(-))}
\]
is the trivial map if $s>t$ and so there is a commutative diagram
\[
\xymatrix{
K_n^R\ar[r]^-{e}\ar[d] & \Coalg(C(-)^{\otimes n},T(-))\ar[d]\\
\frac{K_n^R}{\Gamma^{t+1}K_n^R}\ar[r]^-{e} & \Coalg(\sk_t,T(-))}
\]
for each $t\geq 1$.

Let $E_t$ be the kernel of the map
$\Coalg(\sk_t,T(-))\lra\Coalg(\sk_{t-1},T(-)).$ Then there is a
commutative diagram
\[
\xymatrix{ \Gamma^t/\Gamma^{t+1}\ar@{^(->}[r]\ar[d]^{e} &
K^R_n/\Gamma^{t+1}\ar[d]^{e}\ar@{->>}[r] &
K_n^R/\Gamma^t,\ar[d]^{e}\\
E_t\ar@{^(->}[r]&\Coalg(\sk_t,T(-))\ar[r] &
\Coalg(\sk_{t-1},T(-))}
\]
where $\Gamma^t=\Gamma^tK_n^R.$ Notice that
\[
\Gamma^{n+1}=E_{n+1}=\{1\} \quad \text{ and }\quad
K^R_n/\Gamma^1=\Coalg(\sk_0,T(-))=\{1\}.
\]
It suffices to show that $e\co \Gamma^t/\Gamma^{t+1}\lra E_t$
is an epimorphism for each $t\geq1$.

Notice that $f\in E_t\subseteq\Hom_R(\sk_t,T(-))$ if and only if
$f\co\sk_tC(V)^{\otimes n}\lra T(V)$ is a functorial $R$--linear
map such that
\begin{enumerate}
\item $f|_{(C(V)^{\otimes n})_{j}}=0\co (C(V)^{\otimes
n})_{j}\lra T(V)_j$ for $0<j<2t$; 
\item $f((C(V)^{\otimes
n})_{2t})\subseteq P(T(V)_{2t})=P_t(T(V))$.
\end{enumerate}
Now let $f\in E_t$ and let $\tilde f\co C(V)^{\otimes n}\lra
T(V)$ be the functorial $R$--linear map such that $\tilde
f|_{\sk_t}=f$ and $\tilde f|_{(C(V)^{\otimes n})_j}=0$ for $j>2t$.
Let $1\leq i_1<i_2<\cdots<i_t\leq n$ be a sequence of length $t$.
Let
\[
\lambda_I\co C(V)^{\otimes t}=R\otimes\cdots
R\otimes\stackrel{(i_1)}{C(V)}\otimes R\cdots\otimes R\otimes
\stackrel{(i_t)}{C(V)}\otimes R\cdots\otimes R\lra C(V)^{\otimes
n}
\]
be the inclusion of $i_1,i_2,\cdots,i_t$ coordinates and let
\[
\pi_I=\epsilon\otimes\cdots\otimes\stackrel{(i_1)}{\id_{C(V)}}\otimes\cdots
\otimes\stackrel{(i_t)}{\id_{C(V)}}\otimes\cdots\otimes\epsilon\co
C(V)^{\otimes n}\lra C(V)^{\otimes t}
\]
be the projection of $i_1,i_2,\cdots,i_t$ coordinates. Let
$g_I\co C(V)^{\otimes t}\lra T(V)$ be the composite
\[\xymatrix{
C(V)^{\otimes t}\ar[r]^{\lambda_I}&C(V)^{\otimes n}\ar[r]^{\tilde
f}&T(V)}
\]
and let $\tilde f_I=g_I\circ \pi_I\co C(V)^{\otimes n}\lra
T(V)$. Notice that
\[
\tilde f_I((C(V)^{\otimes n})_{2t})=f_I((C(V)^{\otimes
n})_{2t})\subseteq P(T(V)).
\]
Thus $g_I((C(V)^{\otimes t})_{2t})=g_I(V^{\otimes t})\subseteq
P(T(V))$. Let $\bar V$, $\gamma^R_t$ and $\Lie^R(t)$ be defined as
in \fullref{lemma2.11}. By Lemmas~\ref{lemma5.2}
and~\ref{lemma2.11}, we have
\[
g_I(\gamma^R_t)\subseteq \gamma^R_t\cap P_t(\bar V)=\Lie^R(n).
\]
Thus
\[
g_I(x_1x_2\cdots x_t)=\sum_{\tau\in
S_{n-1}}r_{\tau}[[x_1,x_{\tau(2)},\cdots,x_{\tau(t)}]
\]
for some $r_{\tau}\in R$ for each $\tau\in \Sigma_{n-1}$, where
$\Sigma_{n-1}$ acts on $\{2,3,\cdots,n\}$. Notice that $g_I$ is a
functorial map. Thus
\[
g_I(a_1a_2\cdots a_t)=\sum_{\tau\in
S_{n-1}}r_{\tau}[[a_1,a_{\tau(2)},\cdots,a_{\tau(t)}]
\]
for any free $R$--module $V$ and any $a_j\in V$. Let
\[
w_I=\prod_{\tau\in
S_{n-1}}(1+r_{\tau}[[y_{i_1},y_{i_{\tau(2)}},\cdots,y_{i_{\tau(t)}}])
\in A^R(y_1,y_2,\cdots,y_n),
\]
where the product is taken in an arbitrary order. Then it is
routine to check that
\[
\theta_n(w_I)|_{\sk_t}=g_I\circ \pi_I|_{\sk_t}=\tilde
f_I|_{\sk_t}\co\sk_t C(V)^{\otimes n}\lra T(V).
\]

By \fullref{lemma2.10}, $w_I\in K^R_n$ and so
\[
\tilde f_I|_{\sk_t}\in\Im(e\co\Gamma^t/\Gamma^{t+1}\lra E_t).
\]

Notice that
\[
(C(V)^{\otimes n})_{2t}=\bigoplus_{1\leq i_1<i_2<\cdots<i_t\leq
n}V^{\otimes t}=\bigoplus_{n\choose t}V^{\otimes t}.
\]

Thus
\[
f=\tilde f|_{\sk_t}=\prod_I\tilde f_I|_{\sk_t}
\]
is in the convolution algebra $\Hom_R(\sk_tC(-)^{\otimes
n},T(-))$, where $I=(i_1,i_2\cdots,i_t)$ runs over all sequences
$1\leq i_1<i_2<\cdots<i_t\leq n$ of length $t$ and the product is
taken in an arbitrary order. Hence
\[
f\in\Im(e\co\Gamma^t/\Gamma^{t+1}\lra E_t)
\]
and the assertion follows.
\end{proof}
\medskip

The James filtration of the tensor algebra $T(V)$
\[
J_0(V)\subseteq J_1(V)\subseteq \cdots \subseteq J_n(V)\subseteq
\cdots\subseteq T(V)
\]
induces a cofiltration of the progroup $\Coalg(T(-),T(-))$
\[
\xymatrix{ \cdots\ar@{->>}[r] & \Coalg(J_n(-),T(-))\ar@{->>}[r] &
\cdots \ar@{->>}[r] & \Coalg(J_0(-),T(-)),}
\]
where the group $\Coalg(T(-),T(-))$ can be identify with the
inverse limit
\[
\Coalg(T(-),T(-))=\lim_n \Coalg(J_n(-),T(-)).
\]
Selick and Wu~\cite{SW} proved that there is an isomorphisms of
groups
\[
\xymatrix{\Ker\left[\Coalg(J_n(-),T(-))\right.\ar@{->>}[r]&\left.
\Coalg(J_{n-1}(-),T(-))\right]\ar[r]^-{\cong}&\Lie^{R\otimes\Z/p}(n).}
\]

\begin{defn}
Define the group $\bar L^R_n$ to be the equalizer of the
projection maps $\xymatrix{\pi_j\co G(A^R(y_1,\ldots
y_n))\ar@{->>}[r] & G(A^R(y_1,\ldots y_{n-1}))}$.
\end{defn}
\medskip

\begin{lem}
The group $\bar L^R_n$ is isomorphic to the combinatorial group
$\H_n^R$.
\end{lem}

\begin{proof}
Recall that the isomorphism $K_n\cong G(A^R(y_1,\ldots y_n))$
preserves the projection maps, that is,
$\exp\circ\pi_i=\pi_i\circ\exp$. Therefore the assertion of the
lemma follows immediately as both groups are given by the
equalizer of the same maps.
\end{proof}
\medskip

\begin{thm}\label{theorem3.11}
There is an isomorphism of groups
\[
e\co \H^R_n\lra \Coalg(J_n(-),T(-))
\]
for $1\leq n\leq\infty$.
\end{thm}

\begin{proof}
Let $\theta_n\co L^R_n\lra\Hom_R(J_n(-),T(-))$ be
homomorphism~\eqref{thetaonL}. As the following diagram
\[
\xymatrix{
\H^R_n \ar@{^(->}[r] \ar[d] & L^R_n\ar[d]\\
K_n^R \ar@{^(->}[r]\ar[d]^{p_j} &
A^R(y_1,\cdots,y_n)\ar[d]^{\pi_j}\\
K^R_{n-1}\ar@{^(->}[r] & A^R(y_1,\cdots,y_{n-1}),}
\]
commutes for every $1\leq j\leq n$ and the bottom two rows are
given by inclusions, there is an inclusion $\H^R_n\lra L^R_n$
between the equalizers. Notice that $\theta_n(\H^R_n)\subseteq
\Coalg(J_n(-),T(-))$ and so
\[
e=\theta_n|_{\H^R_n}\co \H_n^R\lra\Coalg(J_n(-),T(-))
\]
is a well-defined monomorphism of groups. We need to show that $e$
is an epimorphism.

Recall that there is an isomorphisms of groups
\[
\xymatrix{ \Ker\big[\Coalg(J_n(-),T(-))\ar@{->>}[r] &
\Coalg(J_{n-1}(-),T(-))\big]\ar[r]^-{\cong} &
\Lie^{R\otimes\Z/p}(n).}
\]
and that the kernel $\Lambda (n)$ of the homomorphism
$d_n\co\H^R_n\lra \H^R_{n-1}$ is by
\fullref{kernel_projectionH_n} isomorphic to $\Lie^R(n)$ as a
group. Thus there is a commutative diagram
\[
\xymatrix{
\Lie^R(n) \ar[r]\ar[d] & \Lie^{R\otimes\Z/p}(n)\ar[d]\\
\H^R_n\ar[r]^-{e}\ar@{->>}[d]^{d_n} &
\Coalg(J_n(-),T(-))\ar@{->>}[d]\\
\H^R_{n-1} \ar[r]^-{e} & \Coalg(J_{n-1}(-),T(-)).}
\]
for $1\leq n\leq\infty$. Therefore, we have that the map $e\co
\Lambda (n)\cong\Lie^R(n)\lra\Lie^{R\otimes\Z/p}(n)$ is onto. Now
by induction, the homomorphism
\[
e\co \H^R_n\lra \Coalg(J_n(-),T(-))
\]
is onto for $0\leq n\leq\infty$. The assertion follows.
\end{proof}
\medskip

\begin{corollary}
There is an isomorphism of groups
\[
e\co \H^R_{\infty}\lra \Coalg(T(-),T(-)).
\]
\end{corollary}
\medskip

Recall that by \fullref{theoremHisomoalg}, the combinatorial
algebra $L^R_n$ is isomorphic to\break $\Hom_R(J_n(-), T(-))$. Therefore
we have the following statement.

\begin{corollary}
There is a faithful representation
\[
e\co \H_n^R\lra L^R_n\quad\text{for $1\leq n\leq\infty$.}
\]
\end{corollary}
\medskip

\begin{thm}\label{theoremHlisogrp}
There is an isomorphism of groups
\[
e\co ^R\H^{(l)}_n\lra \Coalg(J_n(-^{\otimes l}),T(-))
\quad\text{for $1\leq n\leq\infty$.}
\]
\end{thm}
\begin{proof}
The homomorphism $e\co ^R\H^{(l)}_n\lra \Coalg(J_n(-^{\otimes}
\Coalg(J_n(-^{\otimes l}),T(-))$ is the homomorphism of the group
isomorphisms
\[
e\co K^R_n\lra \Coalg(C(-)^{\otimes n}, T(-)),
\]
and thus it is a group isomorphism.
\end{proof}
Recall that by \fullref{thm:Hlkisomoalg} the combinatorial
algebra $^RL^{(l)}_n$ is isomorphic to\break
$\Hom_R(C(-)^{\otimes n},
T(-))$; hence we have the following corollary.
\begin{corollary}
There is a faithful representation
\[
e\co ^R\H^{(l)}_n\lra\,^RL^{(l)}_n
\quad\text{for $1\leq n\leq\infty$.}
\]
\end{corollary}

\section{A representation for the group $K_n(k)$}\label{sec:representationforKn}
In this section, we give a representation for the group $K_n(k)$.

The algebra $A(y_1,\cdots,y_n)$ gives a representation for the
group $K_n=K_n(1)$ (see \fullref{sec:cohen_groups}). We will
use this algebra to give a representation for $K_n(k)$ for general
$k$, where the representation map into $A(y_1,\cdots,y_n)$ depends
on $k$.

Recall that the algebra
$A(y_1,\negthinspace\cdots\negthinspace,y_n)$  is a graded algebra
over $\Z$ and $A(y_1,\negthinspace\cdots\negthinspace,y_n)_t$ is a
free $\Z$--module with basis $ y_{i_{\sigma(1)}}\cdots
y_{i_{\sigma(t)}} $, where $\sigma\in \Sigma_t$ acts on
$\{1,2,\cdots,t\}$ and $1\leq i_1<i_2<\cdots<i_t\leq n$.

The (ungraded) commutator in the algebra $A(y_1,\cdots,y_n)$ is
defined by
\[
[x,y]=xy-yx
\]
for $x,y\in A(y_1,\cdots,y_n)$ and the iterated commutator is
defined inductively by
\[
[[x_1,\cdots,x_t]=[([[x_1,\cdots,x_{t-1}]),x_t].
\]
Now consider the homogeneous monomial $y_{i_1}\cdots y_{i_k}\in
A(y_1,\cdots,y_n)$. Notice that $(y_{i_1}\cdots y_{i_k})^2=0$.
Thus
\[
(1+y_{i_1}\cdots y_{i_k})(1-y_{i_1}\cdots y_{i_k})=1
\]
and so $1+y_{i_1}\cdots y_{i_k}$ is a unit element in
$A(y_1,\cdots,y_n)$. Let $F(k)$ be the free group generated by the
words $\{x_{i_1}|\cdots|x_{i_k}\}$ with $1\leq i_s\leq n$ for
$1\leq s\leq k$ and let
\[
\xymatrix{ e\co F(k)\ar[r] & A(y_1,\cdots,y_n),\mbox{ given by
}\ e(\{x_{i_1}|\cdots|x_{i_k}\})=(1+y_{i_1}\cdots y_{i_k}),}
\]
be a homomorphism from $F(k)$ into the group of units in
$A(y_1,\cdots,y_n)$. Let $I=(i_1,\cdots,i_k)$ be a sequence of
integers. Write $x_I$ for $\{x_{i_1}|\cdots|x_{i_k}\}$ and $y_I$
for the monomial $y_{i_1}\cdots y_{i_k}$.

\begin{thm}\label{representation1}
The representation map $e\co F(k)\lra A(y_1,\cdots,y_n)$
factors through the quotient group $K_n(k)$. Furthermore, the
resulting representation map
\[
e\co K_n(k)\lra A(y_1,\cdots,y_n)
\]
is faithful.
\end{thm}

The proof requires the following lemma.

\begin{lem}\label{lemma3.4}\
\begin{enumerate}
\item Let $w=x^{n_1}_{I_1}\cdots x^{n_t}_{I_t}$ be a word in
$F(k)$. Then
 \[
 e(w)=(1+n_1y_{I_1})\cdots (1+n_t y_{I_t}).
 \]
\item $e([[x_{I_1}^{n_1},\cdots,x_{I_t}^{n_t}])=1+n_1\cdots
n_t[[y_{I_1},\cdots,y_{I_t}]$, where the first bracket $[[\cdots]$
is the commutator in the group $F(k)$ and the second bracket
$[[\cdots]$ is the commutator in the algebra $A(y_1,\cdots,y_n)$.
\end{enumerate}
\end{lem}
\begin{proof}
Notice that $(1+ny_I)(1-ny_I)=1-n^2(y_I)^2=1$ for $n\in\Z$. Thus
\[
e(x^{-n}_I)=1-ny_I
\]
and so assertion (1) follows.

The proof of assertion (2) is given by induction on the length
$t$. Let
\begin{align*}
w&=[[x_{I_1}^{n_1},\cdots,x_{I_t}^{n_t}]
\\
\tag*{\hbox{\rlap{with $t\geq 2$ and let}}}
w'&=[[x_{I_1}^{n_1},\cdots,x_{I_{t-1}}^{n_{t-1}}]
\end{align*}
Notice that $w=(w')^{-1}x^{-n_t}_{I_t}w'x^{n_t}_{I_t}$. By
induction, one has
\begin{align*}
e(w')&=1+n_1\cdots n_{t-1}[[y_{I_1},\cdots,y_{I_{t-1}}]
\\
\tag*{\hbox{and so}}
e((w')^{-1})&=1-n_1\cdots n_{t-1}[[y_{I_1},\cdots,y_{I_{t-1}}]
\end{align*}
because $[[y_{I_1},\cdots,y_{I_{t-1}}]^2=0$ in the algebra
$A(y_1,\cdots,y_n)$. Thus
\[
\begin{array}{l}
e(w)=\\\hspace*{-0.3cm}(1\negthinspace-n_1\negthinspace\cdots\negthinspace
n_{t-1}[[y_{I_1},\negthinspace\cdots\negthinspace,y_{I_{t-1}}])(1\negthinspace-n_t
y_{I_t})(1\negthinspace+n_1\cdots
n_{t-1}[[y_{I_1},\negthinspace\cdots\negthinspace,y_{I_{t-1}}])(1+n_t
y_{I_t}\negthinspace).
\end{array}
\]
It is routine to check that the right side in the above equation
is equal to $$1+n_1\cdots n_t [[y_{I_1},\cdots, y_{I_t}]$$ in the
algebra $A(y_1,\cdots,y_n)$. Assertion (2) follows.
\end{proof}

\begin{proof}[Proof of \fullref{representation1}]
Let $I=(i_1,\cdots,i_k)$. Then $e(x_I)=1+y_{i_1}\cdots y_{i_k}$.
If $i_p=i_q$ for some $p<q$, then $y_{i_1}\cdots y_{i_k}=0$ and so
the representation map $e$ preserves relation~(1) in the
definition of $K_n(k)$.

Let $I_j=(i_{j1},\cdots,i_{jk})$ with $1\leq j\leq t$. Then
\[
e([[x_{I_1},\cdots, x_{I_t}])=1+[[y_{I_1},\cdots,y_{I_t}].
\]
\[
[[y_{I_1},\cdots,y_{I_t}]=\sum_{\sigma}\pm y_{I_{\sigma(1)}}\cdots
y_{I_{\sigma(t)}},\leqno{\hbox{Notice that}}
\]
where $\sigma$ runs over certain elements in $\Sigma_t$.

Suppose that $i_{pq}=i_{p'q'}$ for some $(p,q)\not=(p',q')$. Then
\[
[[y_{I_1},\cdots,y_{I_t}]=0
\]
and so $e([[x_{I_1},\cdots, x_{I_t}])=1$. That is, the
representation map $e$ preserves relation~(2) in the definition of
$K_n(k)$.

Now it is easy to check that
\begin{align*}
(1+[[y_{I_1},\cdots,y_{I_t}])^n&=1+n[[y_{I_1},\cdots,y_{I_t}].
\\
\tag*{\hbox{Thus}}
e([[x_{I_1}^{n_1},\cdots,x_{I_t}^{n_t}])&=1+n_1\cdots
n_t[[y_{I_1},\cdots,y_{I_t}]\\
&=(1+[[y_{I_1},\cdots,y_{I_t}])^{n_1\cdots
n_t}=e([[x_{I_1},\cdots,x_{I_t}]^{n_1\cdots n_t})
\end{align*}
and so the representation map $e$ preserves relation (3) in the
definition of $K_n(k)$. Thus the representation map $e\co
F(k)\lra A(y_1,\cdots, y_n)$ factors through $K_n(k)$, which is
the first statement in the theorem.

Now we need to show that the representation map $e\co
K_n(k)\lra A(y_1,\cdots,y_n)$ is faithful. Let $A(k)$ be the
graded subalgebra of $A(y_1,\cdots,y_n)$ generated by all the
elements $\alpha\in A(y_1,\cdots,y_n)_k$ of dimension $k$. Notice
that the image of the representation map $e$ is contained in the
subalgebra $A(k)$. That is,
\[
e\co K_n(k)\lra A(k)
\]
is a representation map. Let $I(A(k))$ be the augmentation ideal
of the algebra $A(k)$ and let $\theta \co K_n(k)\lra I(A(k))$
be the function $\theta(w)=e(w)-1$ for $w\in K_n(k)$. Let
$\{\Gamma^r(G)\}$ be the lower central series of a group $G$ and
let $I^t(A)=(I(A))^t$ be the $t$--fold product of the augmentation
ideal $I(A)$ of an augmented algebra $A$. To finish the proof the
following two lemmas are required.

\begin{lem}\label{representation1a}
Let $G$ be a group and let $e\co G\lra A$ be a representation
of $G$ into an augmented algebra $A$ over $\Z$ such that
$\epsilon(e(g))=1$ for $g\in G$, where $\epsilon\co A\lra \Z$
is the augmentation map. Let $\theta \co G\to I(A)$ be the
function $\theta(g)=g-1$ for all $g\in G$. Then
\begin{enumerate}
\item $\theta(\Gamma^t(G))\subseteq I^t(A)$ for each $t\geq1$.
\item $\theta$ induces a homomorphism of abelian groups
 \[
 \bar\theta\co \Gamma^t(G)/\Gamma^{t+1}(G)\lra I^t(A)/I^{t+1}(A)
 \]
 such that the diagram
\[
\xymatrix{
\Gamma^t(G)\ar[r]^-{\theta}\ar[d] & I^t(A)\ar[d]\\
\Gamma^t(G)/\Gamma^{t+1}(G)\ar[r]^-{\bar\theta} &
I^t(A)/I^{t+1}(A)}
\]
commutes. \item The map
\[
\xymatrix{
\bar\theta\co\sum_{t\geq1}\Gamma^{t}(G)/\Gamma^{t+1}(G)\ar[r]&
\sum_{t\geq0}I^t(A)/I^{t+1}(A)}
\]
is a representation map of the  Lie algebra
$\sum_{t\geq1}\Gamma^{t}(G)/\Gamma^{t+1}(G)$ into the algebra
$\sum_{t\geq0}I^t(A)/I^{t+1}(A)$.
\end{enumerate}
\end{lem}
\begin{lem}\label{representation1b}
The map
\[
\bar\theta \co\Gamma^t(K_n(k))/\Gamma^{t+1}(K_n(k))\lra
I^t(A(k))/I^{t+1}(A(k))
\]
is a monomorphism for each $t\geq1$.
\end{lem}

\noindent\textbf{Continuation of the Proof of
\fullref{representation1}} (given \fullref{representation1a}
and \fullref{representation1b})\\ By relation (2) in the
definition of the group $K_n(k)$, the group $K_n(k)$ is nilpotent
and so $\Gamma^s(K_n(k))$ is trivial for some~$s$. Suppose that
$e\co \Gamma^{t+1}(K_n(k))\lra A(k)$ is faithful. That is,
$w=1$ if $e(w)=1$ and $w\in \Gamma^{t+1}(K_n(k))$. Let $w\in
\Gamma^t(K_n(k))$ be such that $e(w)=1$. Then $\theta(w)=0$ and so
$\bar\theta(\bar{w})=0$ where $\bar{w}$ is the image of $w$
in~$\Gamma^t\bigl(K_n(k)\bigr)/\Gamma^{t+1}\bigl(K_n(k)\bigr)$. By
\fullref{representation1b}, $\bar{w}$ is zero in the abelian
group $\Gamma^t(K_n(k))/\Gamma^{t+1}(K_n(k))$ and so $w\in
\Gamma^{t+1}(K_n(k))$. By induction, $w=1$ and the assertion
follows.
\end{proof}

\begin{proof}[Proof of \fullref{representation1a}]
Let $\tilde e\co \Z(G)\lra A$ be the map of algebras such that
$\tilde e(g) =$\break $e(g)$, where $\Z(G)$ is the group ring of $G$ over
$\Z$. Notice that $\epsilon(e(g))=1$ for $g\in G$. The map $\tilde
e\co \Z(G)\lra A$ is a morphism of augmented algebras. Let
$\theta'\co G\lra I(\Z(G))$ be the function $\theta'(g)=g-1$
for $g\in G$. Then $\theta=\tilde e\circ\theta'\co G\lra A$.
Notice that the assertions hold for the function
 $\theta'\co G\lra \Z(G)$, (see, for example, \cite[Proposition~10.5,
 page~187]{Curtis}),  and $\tilde e$ is a morphism of augmented algebras.
 The assertion follows.
\end{proof}
We need two more lemmas to prove \fullref{representation1b}.
\begin{lem}{\rm\cite{Liealgebra}}\label{commutator}\qua
Let $A$ be an (ungraded) algebra over $\Z$ and let $x_1,\cdots,
x_n$ be elements in $A$. Then
\[
[[x_1,\cdots,x_n]=\sum_{0\leq p\leq n-1}(-1)^px_{i_p}\cdots
x_{i_1} x_{j_1}\cdots x_{j_{n-p}},
\]
where the sum is taken over $0\leq p\leq n-1$ and all shuffles
$(i_1,\cdots, i_p;j_1,\negthinspace\cdots\negthinspace, j_{n-p})$
such that:
\begin{enumerate}
\item $1\leq i_1<i_2<\cdots<i_p\leq n$; \item $1\leq
j_1<j_2<\cdots<j_{n-p}\leq n$; \item $\{i_1,\cdots,i_p\}\cap
\{j_1,\cdots,j_{n-p}\}=\emptyset$.
\end{enumerate}
\end{lem}
\begin{proof}
The proof is given by induction on the length $n$. The assertion
holds trivially for the case $n=2$. Suppose that the assertion
holds for $n-1$. Then
\[
[[x_1,\cdots,x_n]=[[x_1,\cdots,x_{n-1}]x_n-x_n[[x_1,\cdots,x_{n-1}]
\]
\[
=\negthinspace\sum_{0\leq p\leq
n-2}\hspace{-.2cm}(-1)^px_{i_p}\negthinspace\cdots x_{i_1}
x_{j_1}\cdots x_{j_{n-p-1}}x_n\, -\!\!\!\sum_{0\leq p\leq
n-2} \hspace{-.2cm}(-1)^px_nx_{i_p}\cdots x_{i_1} x_{j_1}\cdots
x_{j_{n-p-1}}
\]
\[ =\hspace*{-0.35cm}\negthinspace\sum_{0\leq p\leq
n-2;i_p<n}\hspace{-.3cm}\negthinspace(-1)^px_{i_p}\negthinspace\cdots
x_{i_1} x_{j_1}\cdots x_{j_{n-p-1}}x_n\, +\hspace{-.2cm}\!\!\!\sum_{1\leq
p\leq
n-1;i_p=n}\hspace{-.2cm}\negthinspace\negthinspace(-1)^{p}x_{i_p}\cdots
x_{i_1}x_{j_1}\cdots x_{j_{n-p}}.
\]
Notice that $x_{i_p}\cdots x_{i_1}x_{j_1}\cdots
x_{j_{n-p}}=x_1\cdots x_n$ or $x_n\cdots x_1$ if $p=0$ or $n-1$,
respectively. The assertion follows.
\end{proof}
A sequence of integers $I=(i_1,\cdots,i_s)$ is called {\it
$n$--admissible} if $i_p\not=i_q$ for all $p<q$ and $1\leq i_j\leq
n$ for $1\leq j\leq s$. The length $s$ of the sequence
$I=(i_1,\cdots,i_s)$ will be denoted by $l(I)$. Let
$I_j=(i_{j1},\cdots, i_{jk_j})$ be sequences of integers with
${1\leq j\leq t}$. The sequence (of sequences) $(I_1,\cdots,I_t)$
is called {\it $n$--admissible} if the sequence $(i_{11},\cdots,
i_{1k_1},\cdots,i_{t1},\cdots,i_{tk_t})$ is $n$--admissible. Now
consider the graded algebra $A(k)$. Let $T(y_I)$ be the tensor
algebra over $\Z$ generated by the words $y_I$, where
$I=(i_1,\cdots,i_k)$ is a sequence such that $1\leq i_j\leq n$ for
$1\leq j\leq k$. Let $q\co T(y_I)\lra A(k)$ be the quotient map
of algebras, where $q(y_I)=y_I$. Let $L(y_I)\subseteq T(y_I)$ be
the associated Lie algebra of $T(y_I)$. Let $A_t(k)$ denote the
image $q(T_t(y_I))$, let $\bar L(A(k))$ denote $q(L(y_I))$ and let
$\bar L_t(A(k))$ denote $q(L_t(y_I))$ where $L_t(y_I)$ is spanned
by the Lie elements of tensor length $t$ in the Lie algebra
$L(y_I)$. Let $\calS_k$ be the set of the $n$--admissible sequences
of length $k$ with some chosen order~$<$.
\eject
\begin{lem}\label{abasis}\
\begin{enumerate}
\item A basis for the free $\Z$--module
$A_t(k)=A(y_1,\cdots,y_n)_{kt}$ is given by the monomials $
y_{I_1}\cdots y_{I_t}, $ where $I_j\in \calS_k$ for $1\leq j\leq
t$ and the sequence $(I_1,\cdots, I_t)$ is admissible. \item A
basis for the free $\Z$--module $\bar L_t(A(k))$ is given by the
Lie elements
\[
[[y_{I_1},y_{I_{\sigma(2)}},\cdots,y_{I_{\sigma(t)}}],
\]
where $I_j\in\calS_k$ for $1\leq j\leq t$ such that
$I_1<I_2<\cdots<I_t$ and where $(I_1,\cdots,I_t)$ is admissible
and $\sigma\in\Sigma_{t-1}$ acts on $\{2,\cdots,t\}$.
\end{enumerate}
\end{lem}
\medskip

\begin{proof}
It is easy to check that the set of the monomials
\[
\{y_{I_1}\cdots y_{I_t}\mid I_j\in \calS_k \mbox{ for } 1\leq
j\leq t {\rm\  and\ } (I_1,\cdots,I_t) \mbox{ is admissible}\}
\]
is the same as that of monomials
\[
\{y_{i_1}\cdots y_{i_{kt}}\mid (i_1,\cdots, i_{kt}) \mbox{ is
admissible}\}.
\]

Now the latter is a basis for $A(y_1,\cdots,y_n)_{kt}=A_t(k)$.
Assertion (1) follows.
\medskip

For assertion (2), it suffices to show that the Lie elements given
in assertion (2) are linearly independent because every element in
$\bar L_t(A(k))$ is a linear combination of such elements.
\medskip

Let $I_j, J_j\in\calS_k$ for $1\leq j\leq t$ such that:
\begin{enumerate}
\item $I_1<\cdots<I_t$, $J_1<\cdots<J_t$; \item both
$(I_1,\cdots,I_t)$ and $(J_1,\cdots,J_t)$ are admissible.
\end{enumerate}
\medskip

Let $\sigma,\tau\in\Sigma_{t-1}$. By assertion (1),
$\Hom_{\Z}(A_t(k),\Z)$ has a dual basis. Consider
\[
\la[[y_{J_1},y_{J_{\tau(2)}},\cdots,y_{J_{\tau(t)}}],
(y_{I_1}y_{I_{\sigma(2)}}\cdots y_{I_{\sigma(t)}})^*\ra.
\]

If the set $\{J_1,\cdots,J_t\}\not=\{I_1,\cdots,I_t\}$, then it is
easy to check that
\[
\la[[y_{J_1},y_{J_{\tau(2)}},\cdots,y_{J_{\tau(t)}}],
(y_{I_1}y_{I_{\sigma(2)}}\cdots y_{I_{\sigma(t)}})^*\ra=0.
\]

If $\{J_1,\cdots,J_t\}=\{I_1,\cdots,I_t\}$ as a set, then
$J_j=I_j$ for all $1\leq j\leq t$. By \fullref{commutator},
\[
[[y_{J_1},y_{J_{\tau(2)}},\cdots,y_{J_{\tau(t)}}]=y_{J_1}y_{J_{\tau(2)}}\cdots
y_{J_{\tau(t)}}+\sum\pm y_{K_1}y_{K_2}\cdots y_{K_t},
\]
where $\sum\pm y_{K_1}y_{K_2}\cdots y_{K_t}$ is the sum of the
other terms in \fullref{commutator}. Notice that
$y_{K_1}\not=y_{J_1}$ in each term of the sum $\sum
y_{K_1}y_{K_2}\cdots y_{K_t}$. Thus the scalar product
\[
\la[[y_{J_1},y_{J_{\tau(2)}},\negthinspace\cdots\negthinspace,y_{J_{\tau(t)}}],
(y_{I_1}y_{I_{\sigma(2)}}\negthinspace\cdots
y_{I_{\sigma(t)}})^*\ra\negthinspace=\negthinspace\left\{
\begin{array}{ll}
1&\mbox{for }(J_1,\negthinspace\cdots\negthinspace,J_t)=(I_1,\negthinspace\cdots\negthinspace,I_t),\\
 & \mbox{and }\tau=\sigma\\
0&\mbox{otherwise}.\\
\end{array}
\right.
\]
Assertion (2) now follows easily.
\end{proof}
\medskip

\begin{proof}[Proof of \fullref{representation1b}]
By \fullref{lemma3.4}, the $\Z$--map
\[
\bar\theta\co \Gamma^t(K_n(k))/\Gamma^{t+1}(K_n(k))\lra
I^t(A(k))/I^{t+1}(A(k))=A_t(k)
\]
maps onto the $\Z$--submodule $\bar L_t(A(k))$. By relation (2) in
the definition of the group $K_n(k)$, every element in the abelian
group $\Gamma^t(K_n(k))/\Gamma^{t+1}(K_n(k))$ is a linear
combination of the elements
$\overline{[[x_{I_1},x_{I_{\sigma(2)}},\cdots,
x_{I_{\sigma(t)}}]}, $ where $I_j\in\calS_k$ for $1\leq j\leq t$,
$(I_1<\cdots<I_t)$ is admissible and $\sigma\in \Sigma_{t-1}$ acts
on $\{2,\cdots,t\}$. Thus the $\Z$--map $\phi\co \bar
L_t(A(k))\lra \Gamma^t(K_n(k))/\Gamma^{t+1}(K_n(k))$ defined by
\[
\phi([[y_{I_1},y_{I_{\sigma(2)}},\cdots,y_{I_{\sigma(t)}}])=
\overline{[[x_{I_1},x_{I_{\sigma(2)}}\cdots,x_{I_{\sigma(t)}}]}
\]
is an epimorphism. Notice that the composite
\[
\xymatrix{ \bar L_t(A(k))\ar@{->>}[r]^-{\phi} &
\Gamma^t(K_n(k))/\Gamma^{t+1}(K_n(k))\ar@{->>}[r]^-{\bar\theta} &
\bar L_t(A(k))}
\]
is an isomorphism. Thus
\[
\phi\co \bar L_t(A(k))\lra
\Gamma^t(K_n(k))/\Gamma^{t+1}(K_n(k))
\]
is a monomorphism and hence it is an isomorphism. The same
argument proves that
\[
\bar\theta\co \Gamma^t(K_n(k))/\Gamma^{t+1}(K_n(k))\lra \bar
L_t(A(k))
\]
is an isomorphism . The assertion follows.
\end{proof}

Consider the composite
\[
\ddr{\bar e\co
K_n(k)}{e}{A(y_1,\cdots,y_n)}{q}{A^{\Z/p^r}(y_1,\cdots,y_n).}
\]
Notice that, in the algebra $A^{\Z/p^r}(y_1,\cdots,y_n)$, we have
\[
(1+y_{i_1}\cdots y_{i_k})^{p^r}=1+p^r y_{i_1}\cdots y_{i_k}=1.
\]
Thus the representation map $\bar e\co K_n(k)\lra
A^{\Z/p^r}(y_1,\cdots,y_n)$ factors through the quotient group
$K^{\Z/p^r}_n(k)$.

\begin{proposition}\label{representation2}
The representation map $\bar e\co K^{\Z/p^r}_n(k)\lra
A^{\Z/p^r}(y_1,\cdots\negthinspace, y_n)$ is faithful.
\end{proposition}
\begin{proof}
Let $A^{\Z/p^r}(k)=A(k)\otimes Z/p^r$. It suffices to show that
\begin{align*}
\bar\theta\co
\Gamma^t(K^{\Z/p^r}_n(k))/\Gamma^{t+1}&(K^{\Z/p^r}_n(k))\\
&\lra
I^t(A^{\Z/p^r}(k))/I^{t+1}(A^{\Z/p^r}(k))\negthinspace=\negthinspace
A_t(k)\otimes \Z/p^r
\end{align*}
is a monomorphism. Notice that
\[
\bar\theta\co
\Gamma^t(K^{\Z/p^r}_n(k))/\Gamma^{t+1}(K^{\Z/p^r}_n(k))\lra
I^t(A^{\Z/p^r}(k))/I^{t+1}(A^{\Z/p^r}(k))
\]
maps onto the submodule $\bar L_t(A(k))\otimes \Z/p^r$. Observe
that $L_t(A(k))\otimes \Z/p^r$ is a free $\Z/p^r$--module with a
basis $[[y_{I_1}, y_{I_{\sigma(2)}},\cdots,y_{I_{\sigma(t)}}]$,
where $I_j\in\calS_k$ for $1\leq j\leq t$, where
$(I_1<\cdots<I_t)$ is admissible and where $\sigma\in
\Sigma_{t-1}$ acts on $\{2,\cdots,t\}$. By the definition of
$K^{\Z/p^r}_n(k)$, the map
\[
\phi\co L_t(A(k))\otimes \Z/p^r\lra
\Gamma^t(K^{\Z/p^r}_n(k))/\Gamma^{t+1}(K^{\Z/p^r}_n(k))
\]
given by
\[
\phi([[y_{I_1},y_{I_{\sigma(2)}}\cdots,y_{I_{\sigma(t)}}])=
\overline{[[x_{I_1},x_{I_{\sigma(2)}}\cdots,x_{I_{\sigma(t)}}]}
\]
is a well-defined epimorphism. The assertion follows from the fact
that the composite $\bar\theta\circ\phi$ is an isomorphism.
\end{proof}

\begin{thm}\label{thm:Kn(k)isogrp}
There is an isomorphism of groups
\[
e\co K_n^R(k)\lra \Coalg(C(-)^{\otimes n},T(-^{\wedge k}))
\quad\text{for $1\leq n\leq\infty$.}
\]
\end{thm}
\begin{proof}
The proof follows along the lines of the proof of
\fullref{theorem3.11}.
\end{proof}

\begin{thm}\label{theoremHlkisogrp}
There is an isomorphism of groups
\[
e\co ^R\H^{(l),(k)}_n\lra \Coalg(J_n(-^{\wedge l}),T(-^{\wedge
k}))
\quad\text{for $1\leq n\leq\infty$.}
\]
\end{thm}
\begin{proof}
The homomorphism $e\co\negthinspace
^R\H^{(l)(k)}_n\negthinspace\to\negthinspace \Coalg(J_n(-^{\wedge}
\Coalg(J_n(-\negthinspace^{\wedge l}),T(-^{\wedge k}))$ is the by
the group isomorphisms
\[
e\co K^R_n(k)\lra \Coalg(C(-)^{\otimes n}, T(-^{\wedge k})),
\]
 and thus it is a group isomorphism.
\end{proof}
\begin{corollary}
There is a faithful representation
\[
e\co ^R\H^{(l),(k)}_n\lra\,^RL^{(l),(k)}_n
\quad\text{for $1\leq n\leq\infty$.}
\]
\end{corollary}

\bibliographystyle{gtart}
\bibliography{link}

\end{document}